\documentclass[12pt,centertags,oneside]{amsart}
\usepackage{amsmath,amstext,amsthm,amscd,typearea,hyperref}
\usepackage{amssymb}
\usepackage{a4wide}
\usepackage[mathscr]{eucal}
\usepackage{mathrsfs}
\usepackage{typearea}
\usepackage{charter}
\usepackage{pdfsync}
\usepackage[a4paper,width=16.2cm,top=3cm,bottom=3cm]{geometry}

\numberwithin{equation}{section}

\newtheorem{theorem}{Theorem}[section]
\newtheorem{definition}[theorem]{Definition}
\newtheorem{proposition}[theorem]{Proposition}
\newtheorem{corollary}[theorem]{Corollary}
\newtheorem{lemma}[theorem]{Lemma}
\newtheorem{remark}[theorem]{Remark}

\newtheorem{problem}[theorem]{Problem}

\newcommand{\cali}[1]{\mathscr{#1}}

\newcommand{\GL}{{\rm GL}}

\newcommand{\U}{{\rm U}}

\newcommand{\lof}{\mathop{\mathrm{{log^\star}}}\nolimits}
\newcommand{\supp}{{\rm supp}}

\newcommand{\dist}{\mathop{\mathrm{dist}}\nolimits}

\newcommand{\ddc}{dd^c}
\newcommand{\dc}{d^c}

\newcommand{\dbar}{\overline\partial}
\newcommand{\ddbar}{\partial\overline\partial}
\renewcommand{\GL}{{\rm GL}}

\newcommand{\id}{{\rm id}}

\newcommand{\codim}{{\rm codim\ \!}}

\newcommand{\Ac}{\cali{A}}

\newcommand{\Cc}{\cali{C}}
\newcommand{\Dc}{\cali{D}}
\newcommand{\Ec}{\cali{E}}
\newcommand{\Fc}{\cali{F}}

\newcommand{\Lc}{\cali{L}}

\newcommand{\B}{\mathbb{B}}
\newcommand{\C}{\mathbb{C}}
\newcommand{\D}{\mathbb{D}}
\newcommand{\N}{\mathbb{N}}

\newcommand{\K}{\mathbb{K}}
\newcommand{\R}{\mathbb{R}}
\renewcommand\P{\mathbb{P}}
\newcommand{\T}{\mathbb{T}}
\renewcommand{\U}{\mathbb{U}}

\title[Ergodic theory for Riemann  surface  laminations]{Ergodic theory for Riemann  surface  laminations: a  survey}

 \author{Vi{\^e}t-Anh Nguy{\^e}n}
\address{Universit\'e de Lille 1, 
Laboratoire de math\'ematiques Paul Painlev\'e, 
CNRS U.M.R. 8524,  
59655 Villeneuve d'Ascq Cedex, 
France.}
\email{Viet-Anh.Nguyen@math.univ-lille1.fr,
{\tt http://www.math.univ-lille1.fr/$\sim$vnguyen}}
 
\date{April 12, 2018}
\dedicatory{Dedicated to Professor Kang-Tae Kim for his sixtieth birthday}
\begin{document}

\begin{abstract}
  We  survey   some recent developments in the  ergodic  theory for hyperbolic  Riemann  surface laminations.
The emphasis  is  on
 singular  holomorphic  foliations. These results not  only illustrate  the  strong  similarity  between the ergodic  theory of maps
and that  of  Riemann  surface laminations, but  also indicate  the  fundamental  differences between these two theories.
\end{abstract}

\maketitle

\medskip

\noindent
{\bf Classification AMS 2010}: Primary: 37A35, 37A30, 57R30; Secondary: 58J35, 58J65,
60J65.

\medskip


\tableofcontents


 
\section{Introduction} \label{S:introduction}
 
 These notes   are based  on a series of lectures given by the   author  at KIAS  and at  the KSCV Symposium 11 in Gyeongju in 2016.
 The  purpose is  to review  some   developments in the ergodic theory of laminations by  hyperbolic Riemann surfaces. In particular, 
 we   focus on   the ergodic theory of singular holomorphic foliations.  The   emphasis is  on recent results, but we also include some classical ones for  the sake of completeness and historical
 perspective.

 There is  a  well-known connection between   Riemann  surface laminations  and  the dynamics of iterations  of continuous maps.
 In the  meromorphic  context, this becomes  a link between singular holomorphic  foliations by Riemann surfaces  in dimension $k\geq 2$ and 
 the dynamics of iterations of meromorphic maps in dimension  $k-1.$
 On the one side,
  the abstract ergodic theory of  maps
  has reached  maturity  with remarkable achievements like the Oseledec-Pesin theory.
 The  ergodic  theory of   meromorphic maps  is, however,  much less  developed. Indeed, it has only been  studied  intensively  during the last three decades, see  the survey of Dinh-Sibony
 \cite{DinhSibony10}.
 On the other  side, the ergodic  theory for  hyperbolic  Riemann surface laminations, and in particular, for  the subclass consisting of  singular holomorphic foliations by Riemann surfaces,
 is only in the
early stages of development  and faces a range of challenges in finding its  own way.  In this paper we
describe some  recent  approaches to this  new theory.
We hope that  the ideas  reviewed  in this notes  will  be    developed and expanded   in the  future. 
In writing these notes,  we are  inspired  by   the  surveys and lecture notes of Deroin   
 \cite{Deroin13}, 
    Forn\ae ss-Sibony \cite {FornaessSibony08},  Ghys \cite{Ghys}, Hurder \cite{Hurder}, Zakeri  \cite{Zakeri} etc. In particular, we  are largely influenced  by  
the survey of  Forn\ae ss-Sibony \cite {FornaessSibony08} which gives  an  introduction to  harmonic  currents on singular foliations as developed by themselves. Harmonic  currents are  the analog of
invariant measures in  discrete dynamics. Their approach
 opens  new  avenues  in  
studying the  interplay  between geometry, topology and  dynamics   in the theory of hyperbolic Riemann  surface laminations. 
 
  In Section  \ref{S:basics}, we  will recall  basic facts on Riemann  surface laminations (without  and  with  singularities),  singular holomorphic foliations.
 The hyperbolicity and the leafwise Poincar\'e metric will be introduced. As  consequences, we will develop the  heat diffusions and define the notion of harmonic  measures
 for  hyperbolic Riemann  surface laminations. We also  recall from \cite{FornaessSibony08}  the notion of  positive harmonic currents directed by a Riemann surface lamination
 (possibly with singularities),
 and compare it with   the notion of positive  harmonic  currents on complex  manifolds.  We  give  a  short digression to the isolated singularities for singular holomorphic foliations.
  Singular  holomorphic foliations by Riemann surfaces in $\P^k$ $(k>1)$ provide a large family of examples  where all the above notions  apply.
  In the light of   recent results of Jouanolou \cite{Jouanolou}, Lins Neto-Soares \cite{NetoSoares}, Glutsyuk  \cite{Glutsyuk},  Lins Neto \cite{Neto}, Brunella \cite{Brunella},
  and Loray-Rebelo \cite{LorayRebelo},
  we will  describe  the  properties of a generic  holomorphic foliation in $\P^k$  with  a  given degree $d>1.$

  In Section  \ref{S:Poincare}
  we will introduce  a function  $\eta$  which measures the ratio between the ambient metric and the leafwise  Poincar\'e metric of a lamination.
  This  function plays an important role in the study of  laminations by hyperbolic  Riemann  surfaces.
  We  also  introduce  the class of  Brody hyperbolic  laminations.  This class contains not only all  compact laminations by hyperbolic Riemann surfaces,
  it also includes  many  singular  holomorphic  foliations.
  We then  state some recent results on the regularity  of Brody hyperbolic laminations which arise  from  our  joint-works with Dinh and Sibony in \cite{DinhNguyenSibony14a,DinhNguyenSibony14b}.  
  
  In Section  \ref{S:Mass}
  we study  the mas-distribution  for directed  positive harmonic currents in the local and global settings. Applications to the recurrence  phenomenon of a generic leaf will
  be considered. The material for this  section is mainly  taken from  \cite{DinhNguyenSibony12,NguyenVietAnh17c}.
  
  A  fundamental  contribution to the ergodic theory of laminations/foliations has been made by Garnett in \cite{Garnett} where she introduces the notion of directed harmonic measures
  and  considers the diffusions of the  heat equation in the Riemannian context.
  This idea is further  developed by  Candel \cite{Candel2}.
In Section \ref{S:Ergodic_theorems} we introduce the  diffusions  of the  heat equation for laminations  (possibly
with singularities)   with  respect to a  positive harmonic current directed by a lamination.
This approach  allows us in \cite{DinhNguyenSibony12} to extend the classical theory of  Garnett \cite{Garnett} and  Candel \cite{Candel2}
to Riemann surface laminations  with singularities or to foliations with not
necessarily bounded geometry.  We present two kinds of
ergodic theorems for such currents: one associated to the heat
diffusions and one of  geometric nature close to Birkhoff's averaging on orbits of a
dynamical system.

In Section  \ref{S:Entropy}  
we present a notion  of hyperbolic entropy, using hyperbolic time, for laminations by hyperbolic Riemann surfaces. 
When the lamination is compact and transversally smooth, we state some theorems on the finiteness of   the hyperbolic entropy.
A notion of metric entropy is also introduced for directed  positive harmonic measures.  
This  section is  based on our  joint-works with Dinh and Sibony in \cite{DinhNguyenSibony14a,DinhNguyenSibony14b}.

  Section \ref{S:Lyapunov} is devoted to the Lyapunov theory for hyperbolic Riemann surface laminations.
The  central objects of this theory are the  cocycles  which   are modelled  on the  holonomy cocycle of a  foliation.
We state the  Oseledec multiplicative ergodic theorem  for laminations.  Next,  we  apply it to  smooth compact laminations by hyperbolic Riemann surfaces   and to  compact singular 
holomorphic foliations by Riemann  surfaces. After all, we  characterize    geometrically the Lyapunov exponents of a  smooth cocycle with respect to a harmonic measure.
This  section  is a synthesis of  our several  works in \cite{NguyenVietAnh17a,NguyenVietAnh17b,NguyenVietAnh18}.

Several  open  problems develop in the  course of the  exposition.  Finally, since the choice of the material reflects the limited knowledge of the author
on the ergodic theory of hyperbolic Riemann surface  laminations,
we note that  many  topics 
are not included  here.   The  author apologizes in  advance  for  omissions or undue biases, and  will  welcome comments of suggested inclusions.
 
\bigskip
\noindent
{\bf Main notation.} Throughout the paper, $\D$ denotes the unit disc
in $\C$, $r\D$ denotes the disc of center 0 and of radius $r,$ and
$\D_R\subset\D$ is the disc of center $0$ and of radius $R$ with
respect to the Poincar{\'e} metric on $\D$,
i.e. $\D_R=r\D$ with $R:=\log[(1+r)/(1-r)]$. Poincar{\'e} metric on
a hyperbolic Riemann surface, in particular on $\D$ and on the leaves of a hyperbolic Riemann  surface lamination, is
given by a positive $(1,1)$-form that we denote by $g_P$. The associated distance  is denoted by $\dist_P.$ 
Given a  Riemann surface lamination $(X,\Lc),$
a leaf through a point $x\in X$ is often denoted by $L_x.$   
 Recall that $\dc:={i\over2\pi}(\dbar-\partial)$ and $\ddc = {i\over \pi}\ddbar$.

\smallskip

\noindent{\bf Acknowledgement.}
I would like to thank Nessim Sibony and Tien-Cuong Dinh  for their constant support and for their collaborations  which are largely reported in  these notes.
Sincere thanks also go to  Kyewon Koh, Sung Yeon Kim, Kang-Hyurk Lee, Jong-Do Park, Hyeseon Kim and Taeyong Ahn
  for  their  hospitality and
  encouragement.
   I am  also  grateful to  Mihai Paun, Jaigyoung Choe  and  Lee Kang Won for very kind  help.
The paper was partially prepared 
during my visit  at  Vietnam  Institute for Advanced Study in Mathematics (VIASM) and 
at the Center for Mathematical Challenges (CMC) of the Korea Institute for Advanced Study (KIAS).
  I would like to express my gratitude to these organizations for hospitality and  for  financial support.

\section{Basic results} \label{S:basics}


\subsection{Riemann surface laminations}\label{SS:RSL}

Let $X$ be a locally compact space.  A {\it   Riemann surface lamination}     $(X,\Lc)$   is  the  data of  a {\it (lamination)  atlas} $\Lc$ 
of $X$ with (laminated) charts 
$$\Phi_p:\U_p\rightarrow \B_p\times \T_p.$$
Here, $\T_p$ is a locally compact  metric space, $\B_p$ is a domain in $\C$,  $\U_p$ is  an open set in 
$X,$ and  
$\Phi_p$ is  a homeomorphism,  and  all the changes of coordinates $\Phi_p\circ\Phi_q^{-1}$ are of the form
$$x=(y,t)\mapsto x'=(y',t'), \quad y'=\Psi(y,t),\quad t'=\Lambda(t),$$
 where $\Psi,$ $\Lambda$ are continuous  functions and $\Psi$ is  holomorphic in  $y.$

The open set $\U_p$ is called a {\it flow
  box} and the Riemann surface $\Phi_p^{-1}\{t=c\}$ in $\U_p$ with $c\in\T_p$ is a {\it
  plaque}. The property of the above coordinate changes insures that
the plaques in different flow boxes are compatible in the intersection of
the boxes. Two plaques are {\it adjacent} if they have non-empty intersection.
 
A {\it leaf} $L$ is a minimal connected subset of $X$ such
that if $L$ intersects a plaque, it contains that plaque. So a leaf $L$
is a  Riemann surface  immersed in $X$ which is a
union of plaques. For every point $x\in X,$  denote  by  $L_x$   the   leaf passing  through $x.$  
  A subset $M\subset X$  is  called
  {\it leafwise  saturated} if $x\in M$ implies  $L_x\subset M.$

 We say  that a Riemann  surface lamination  $(X,\Lc)$ is  {\it smooth} if
  each map $\Psi$ above is smooth   with
respect to $y,$ and its partial derivatives of any order  with respect to $y$ and $\bar y$ are jointly continuous
with respect  to $(y,t).$
 
We are mostly interested in the case where the $\T_i$ are closed subsets of smooth real manifolds and the functions $\Psi,\Lambda$ are smooth in all variables. 
In this case, we say that the lamination $(X,\Lc)$ is  {\it transversally smooth}.
If, moreover, $X$ is compact, we can embed it in an $\R^N$ in order to use the distance induced by a Riemannian metric on $\R^N$. 

 We say that a transversally smooth Riemann surface  lamination $(X,\Lc)$ is a  {\it smooth foliation} if
 $X$ is a  manifold and all leaves of $\Lc$ are Riemann surfaces immersed in $X.$
 
We  say that  a Riemann  surface lamination $(X,\Lc)$ is a {\it holomorphic  foliation} 
  if $X$ is a complex manifold (of dimension $k$) and
there is an atlas $\Lc$ of $X$ with charts 
$$\Phi_i:\U_i\rightarrow \B_i\times \T_i,$$
where the $\T_i$'s  are open sets of  $\C^{k-1}$ and 
 all above maps  $\Psi,\Lambda $ are   holomorphic. 
 
 Many examples of abstract compact Riemann  surface laminations are constructed in 
\cite{CandelConlon2} and \cite{Ghys}. Suspensions of  group actions give already a large class
of laminations without singularities.

\subsection{Hyperbolicity and leafwise Poincar\'e metric}\label{SS:Poincare}
 Consider now a  Riemann surface lamination $(X,\Lc)$. 
 \begin{definition}\label{D:hyperbolic}\rm 
A leaf $L$  of $(X,\Lc)$ is  said to be  {\it hyperbolic} if
it  is a   hyperbolic  Riemann  surface, i.e., it is  uniformized   by 
$\D.$   $(X,\Lc)$   is  said to be {\it hyperbolic} if  
  its leaves   are all  hyperbolic. 
\end{definition}
 \smallskip

For every $x\in X$  such that $L_x$ is  hyperbolic,  consider a universal covering map
\begin{equation}\label{e:covering_map}
\phi_x:\ \D\rightarrow L_x\qquad\text{such that}\  \phi_x(0)=x.
\end{equation}
 This map is
uniquely defined by $x$ up to a rotation on $\D$. 
Then, by pushing   forward  the Poincar\'e metric $g_P$
on $\D$  
  via $\phi_x,$ we obtain the  so-called {\it Poincar\'e metric} on $L_x$ which depends only on the leaf.  
  The latter metric is given by a positive $(1,1)$-form on $L_x$  that we also denote by $g_P$ for the sake of simplicity.

 \subsection{Heat diffusions and harmonic  measures}
 \label{ss:heat_diffusions}
     
Let  $(X,\Lc)$ be  a hyperbolic  Riemann surface  lamination.  
     The leafwise Poincar\'e metric
$g_P$  induces  the corresponding 
Laplacian $\Delta$  on leaves   (see formula \eqref{e:Laplacian} for $\beta:=g_P$ below).  
 For  every point  $x\in X,$
 consider  the   {\it heat  equation} on $L_x$
 $$
 {\partial p(x,y,t)\over \partial t}=\Delta_y p(x,y,t),\qquad  \lim_{t\to 0+} p(x,y,t)=\delta_x(y),\qquad   y\in L_x,\ t\in \R^+.
 $$
Here   $\delta_x$  denotes  the  Dirac mass at $x,$ $\Delta_y$ denotes the  Laplacian  $\Delta$ with respect to the  variable $y,$
 and  the  limit  is  taken  in the  sense of distribution, that is,
$$
 \lim_{t\to 0+
}\int_{L_x} p(x,y,t) f(y) g_P( y)=f(x)
$$
for  every  smooth function  $f$   compactly supported in $L_x.$   

The smallest positive solution of the  above  equation, denoted  by $p(x,y,t),$ is  called  {\it the heat kernel}. Such    a  solution   exists   because  $(L_x,g_P)$ is
complete and   of bounded  geometry  (see, for example,  \cite{CandelConlon2,Chavel}).  
 The  heat kernel  $p(x,y,t)$  gives  rise to   a one  parameter  family $\{D_t:\ t\geq 0\}$ of  diffusion  operators    defined on bounded measurable functions  on $X$ by
 \begin{equation}\label{e:diffusions}
 D_tf(x):=\int_{L_x} p(x,y,t) f(y) g_P (y),\qquad x\in X.
 \end{equation}
 This  family is  a  semi-group, that is,  
\begin{equation}\label{e:semi_group}
D_0=\id\quad\text{and}\quad   D_t \mathbf{1} =\mathbf{1}\quad\text{and}\quad D_{t+s}=D_t\circ D_s \quad\text{for}\ t,s\geq 0,
\end{equation}
where $\mathbf{1}$ denotes the function which is identically equal to $1.$

We denote by  $\Cc(X,\Lc)$ the  space of all functions  $f$ defined and  compactly supported  on $ X$ which are  leafwise  $\Cc^2$-smooth
 and transversally continuous, that is,  
 for each  laminated   chart
$\Phi_p:\U_p\rightarrow \B_p\times \T_p$ and all $r,s\in\N$ with $r+s\leq 2,$
 the derivatives  ${\partial^{r+s}(f\circ\Phi_p^{-1})\over  \partial y^r\partial\bar{y}^s}$ exist and are jointly continuous in $(y,t).$

 \begin{definition}\label{D:harmonic_measure}\rm 
 Let  $\Delta$  be   the  Laplacian  on $\Delta,$ that is, the  aggregate  of the  leafwise Laplacians  $\{\Delta_x\}_{x\in X}.$
 
A positive  Borel measure  $\mu$ on $X$ is said  to be
{\it quasi-harmonic}  if
$$
\int_X  \Delta u \,d\mu=0
$$ 
 for all  functions  $u\in \Cc(X,\Lc).$
 
 A  quasi-harmonic measure $\mu$ is  said to be {\it harmonic }  if  $\mu$ is finite and 
 $\mu$  is {\it $D_t$-invariant} for all $t\in\R^+,$ i.e, 
$$ \int_X  D_tf d\mu=\int_X fd\mu,  \qquad  f\in \Cc(X,\Lc),\ t\in\R^+.   $$
\end{definition}


 \subsection{Positive harmonic currents on complex manifolds}
 \label{SS:Hamonic_currents}

 Let $M$ be a  complex manifold of dimension $k.$ 
A $(p,p)$-form on  $M$  is {\it
  positive} if it can be written at every point as a combination with
positive coefficients of forms of type
$$i\alpha_1\wedge\overline\alpha_1\wedge\ldots\wedge
i\alpha_p\wedge\overline\alpha_p$$
where the $\alpha_j$ are $(1,0)$-forms. A $(p,p)$-current or a $(p,p)$-form $T$ on $M$ is
{\it weakly positive} if $T\wedge\varphi$ is a positive measure for
any smooth positive $(k-p,k-p)$-form $\varphi$. A $(p,p)$-current $T$
is {\it positive} if $T\wedge\varphi$ is a positive measure for
any smooth weakly positive $(k-p,k-p)$-form $\varphi$. 
If $M$ is given with a Hermitian metric $\beta$ and $T$ is  a positive  $(p,p)$-current on $M,$
$T\wedge \beta^{k-p}$ is a positive measure on
$M$. The mass of  $T\wedge \beta^{k-p}$
on a measurable set $E$ is denoted by $\|T\|_E$ and is called {\it the mass of $T$ on $E$}.
{\it The mass} $\|T\|$ of $T$ is the total mass of  $T\wedge \beta^{k-p}$ on $M.$

 A $(p,p)$-current on  $M$  is {\it
  harmonic} if  $\ddc T=0$ in the  weak sense (namely,  $T(\ddc f)=0$ for all compactly  smooth $(k-p-1,k-p-1)$-forms  $f$ on $M$). 

  In this article,  for  every $r>0$ let $B_r$ denote the ball of
  center $0$ and of radius $r$ in $\C^k$. 
The following local property of positive harmonic
currents  is discovered  by Skoda \cite{Skoda}.

\begin{proposition} {\rm (Skoda \cite{Skoda}).} \label{P:Skoda}
Let $T$ be a positive harmonic current in
  a ball $B_{r_0}$.
Define $\beta:=\ddc\|z\|^2$
  the standard K{\"a}hler form where $z$ is the canonical
  coordinates on $\C^n$.
Then the function $r\mapsto \pi^{-(k-p)} r^{-2(k-p)}\|T\wedge\beta^{k-p}\|_{B_r}$
is increasing on $0<r\leq r_0$. In particular, it is bounded on 
$]0,r_1]$ for any  $0<r_1<r_0$.
\end{proposition}

The limit of the above function when $r\rightarrow 0$ is called {\it
  the Lelong number} of $T$ at $0$. The above proposition shows that Lelong
number always  exists and is finite positive.

The  next  simple  result    allows  for extending  positive harmonic currents    through isolated points.
\begin{proposition} \label{P:extension}{\rm  (Dinh-Nguyen-Sibony \cite[Lemma 2.5]{DinhNguyenSibony12})}
Let $T$ be a positive current of bidimension $(1,1)$ with compact support on a complex manifold $M$. Assume that $\ddc T$ is a negative measure on $M\setminus E$ where $E$ is a finite set. 
Then  $T$ is a positive  harmonic  current on $M.$
\end{proposition}

 \subsection{Directed positive harmonic currents}
 \label{SS:Directed_hamonic_currents}

  Let $(X,\Lc)$ be  a Riemann surface lamination.
  Let $\Cc^{1,1}(X,\Lc)$ denote the  space   of all forms $f$ of bidegree $(1,1)$   defined  on
leaves  of the lamination and  compactly supported  on $X$  such that $f$ is    transversally continuous.
The last continuity condition means that
 for each  laminated   chart
$\Phi_p:\U_p\rightarrow \B_p\times \T_p,$ the form
   $f\circ\Phi_p^{-1}$ is jointly continuous in $(y,t).$ 
  For each chart $\Phi_p:\U_p\rightarrow \B_p\times \T_p,$ the  complex   structure on $\B_p$ induces  a complex  structure on the leaves of $X.$
Therefore,  the operator $d$  and $\dc$  can be  defined  so that they act leafwise on  forms
as in the case of complex manifolds. 
 So we get easily   that
$\ddc:\  \Cc(X,\Lc)\to\Cc^{1,1}(X,\Lc)$. 
   A form $f\in \Cc^{1,1}(X,\Lc)$ is  said to be {\it positive} if its restriction to every plaque
 is  a  positive $(1,1)$-form in the  usual  sense.
 \begin{definition}\label{D:Directed_hamonic_currents} {\rm (Garnett \cite{Garnett}, see also  Sullivan \cite{Sullivan}).}
\rm A {\it directed  current }  on $(X,\Lc)$  (or equivalently,  a {\it current directed  by the lamination} $(X,\Lc)$) is a   linear continuous  form
     on $\Cc^{1,1}(X,\Lc).$ Let $T$ be a directed current. 
    
  $\bullet$  $T$ is  said to be  {\it positive} if  $T(f)\geq 0$ for all positive forms $f\in \Cc^{1,1}(X,\Lc).$
    
$\bullet$ $T$ is  said to be  {\it   harmonic current} $T$  if  $\ddc T=0$ in the  weak sense (namely,  $T(\ddc g)=0$ for all functions  $g\in  \Cc(X,\Lc)$).
\end{definition}
  We have the  following decomposition.
  \begin{proposition}\label{P:decomposition}
   Let $T$ be  a directed positive harmonic current on $(X,\Lc).$  Let $\U\simeq \B\times \T$ be  a flow  box which is relatively compact in $X.$
   Then, there is a positive Radon measure $\nu$ on $\T$ and for $\nu$-almost every $t\in\T,$ there is a positive  harmonic   function $h_t$ on $\B$
   such that  if $K$  is compact in $\B,$  the integral $\int_K \|h_t\|_{L^1(K)} d\nu(t)$ is  finite and
   $$
   T(f)=\int_\T\big( \int_\B  h_t(y)f(y,t)  \big)d\nu(t)
   $$
   for every  form $f\in\Cc^{1,1}(X,\Lc)$ compactly  supported  on $\U.$ 
  \end{proposition}

\subsection{Directed positive harmonic currents vs harmonic  measures}

Recall that a positive finite  measure  $\mu$ on the $\sigma$-algebra of Borel sets in $X$ is  said  to be  {\it ergodic} if for every  leafwise  saturated  Borel measurable set $Z\subset X,$
   $\mu(Z)$ is  equal to either $\mu(X)$ or $0.$   A   directed positive harmonic  current $T$ is said to be {\it extremal}
   if  $T=T_1+T_2$ for   directed positive harmonic  current $T_1,T_2$ implies that $T_1=\lambda T$  for some $ \lambda \in  [0, 1].$  
The  following result relates the notions of harmonic  measures and  directed positive harmonic currents (see  \cite{DinhNguyenSibony12,NguyenVietAnh17b}). 
\begin{theorem}\label{thm_harmonic_currents_vs_measures} 
Let $(X,\Lc)$ be  a hyperbolic Riemann surface  lamination. \\
(i) If $X$ is  compact, then  each quasi-harmonic    measure is  harmonic.
\\ 
(ii)   The map $T\mapsto  \mu=T\wedge g_P$ which is defined on the  convex  cone of all  directed positive harmonic  currents is one-to-one 
and  its image is contained in  the convex  cone of  all quasi-harmonic  measures $\mu$.
If, moreover, $X$ is  compact, then
 this map  is  a bijection   from the  convex  cone of  all  directed positive harmonic  currents $T$  onto   
the convex  cone of all harmonic  measures $\mu$. 
\\
(iii) If $T$ is an  extremal  directed positive harmonic current and $\mu:=T\wedge g_P$ is  finite, then    $\mu$ is  ergodic.
\end{theorem}

\subsection{Riemann surface laminations with singularities, singular holomorphic  foliations and examples} \label{SS:RSLS}

We call {\it Riemann surface lamination with singularities} the data
$(X,\Lc,E)$ where $X$ is a locally compact space, $E$ a closed
subset of $X$ and $(X\setminus E,\Lc)$ is a Riemann surface
lamination. The set $E$ is {\it the singularity set} of the lamination.
In order to simplify the presentation, we will mostly consider the
case where $X$ is a closed 
subset of a
complex manifold $M$ of dimension $k\geq 1$ and $E$ is a locally finite subset
of $X$. 
We assume that $M$ is endowed with a Hermitian metric $g_M$. 
We also assume that the complex structures on
the leaves of the lamination coincide with the ones induced by $M$,
that is, the leaves of $(X\setminus E,\Lc)$ are Riemann surfaces holomorphically immersed in
$M$. 

 We say that $\Fc:=(X,\Lc,E)$ is  a {\it singular  foliation} (resp.   {\it singular holomorphic foliation}) if $X$ is   a  complex manifold and $E\subset X$ is  a closed subset such that
 $\overline{X\setminus E}=X$ and  $(X\setminus E,\Lc)$ is a  foliation (resp. a  holomorphic  foliation). $E$ is  said to be the  {\it set of singularities} of the  foliation $\Fc.$ 
 We  say that $\Fc$ is  compact if  $X$ is compact.
 
 \begin{definition}\rm
  Let $Z=\sum_{j=1}^k F_j(z){\partial\over\partial  z_j}$ be  a holomorphic vector field  defined in a neighborhood $U$ of $0\in\C^k.$ Consider  the holomorphic map
  $F:=(F_1,\ldots,F_k):\  U\to\C^k.$ We say that $Z$ is 
\begin{enumerate}
\item {\it  singular  at $0$} if  $F(0)=0.$
\item {\it generic linear} if
it can be written as 
$$Z(z)=\sum_{j=1}^k \lambda_j z_j {\partial\over \partial z_j}$$
where $\lambda_j$ are non-zero complex numbers.

\item     {\it with non-degenerate singularity at $0$} if $Z$ is  singular at $0$ and the eigenvalues $\lambda_1,\ldots,\lambda_k$ of the Jacobian matrix $DF(0)$
 are all  nonzero.

\item    {\it with hyperbolic  singularity at $0$} if $Z$ is  singular at $0$ and the eigenvalues $\lambda_1,\ldots,\lambda_k$ of the Jacobian matrix $DF(0)$
 satisfy  $\lambda_j\not=0$ and $\lambda_i/\lambda_j\not\in\R$  for all $1\leq i\not=j\leq k.$
\end{enumerate}
The integral curves of $Z$ define a  singular holomorphic foliation on
$U$.
The condition $\lambda_j\not=0$ implies that the foliation   
has an isolated singularity at 0.
\end{definition}

 Let $\Fc=(X,\Lc,E)$ be a singular holomorphic foliation  such that $E$ is an  analytic subset of $X$  with $\codim(E)\geq 2.$ Then $\Fc$ is    given locally  by holomorphic  vector fields and its
 leaves  are locally, integral curves  of these vector fields, and the singularities of $\Fc$  coincide  with the  singular set of these vector fields.
 We say that a
singular point $a\in E $ is {\it linearizable} (resp.  {\it hyperbolic}) if 
there is a local holomorphic coordinate system of $X$ near $a$ on which the
leaves of $\Fc$ are integral curves of a generic linear vector field  (resp.  of a holomorphic vector field admitting $0$ as  a hyperbolic singularity). 
In   dimension 2 (i.e.  $\dim X=2$),  if  $a$ is  a hyperbolic singularity, then there is a local holomorphic coordinates system of $X$ near $a$ on which the
leaves of $\Fc$ are integral curves of a  vector field $Z(z_1,z_2) = \lambda_1 z_1 {\partial\over \partial z_1}
+ \lambda_2 z_2{\partial\over \partial z_2},$ where  $\lambda_1,\lambda_2$ are some nonzero complex numbers with $\lambda_1/\lambda_2\not\in \R.$ 
In particular, $a$ is a linearizable  singularity.
The   analytic curves $\{  z_1=0\} $ and   $\{  z_2=0\}$ are  called  {\it separatrice}  at $a.$

 Now  we discuss  singular holomorphic  foliations  on $\P^k$ with $k\geq 2.$ Let $\pi:\ \C^{k+1}\setminus\{0\}\to\P^k$ denote  the  canonical projection.  Let $\Fc$ be a
 singular holomorphic  foliations  on $\P^k,$ 
 It can be shown that  $\pi^*\Fc$ is  a singular foliation on $\C^{k+1}$ associated to a vector field $Z$  of the  form
 $$
 Z:=\sum_{j=0}^k F_j(z){\partial\over\partial  z_j},
 $$
 where the  $F_j$ are homogeneous  polynomials of degree $d\geq 1.$ We call  $d$ the {\it degree} of the  foliation. A point $x\in \P^k$
 is a   singularity of $\Fc$ if $F(x)$ is colinear with $x,$ i.e., if $x$ is either an indeterminacy point or a fixed  point  of
 $f=[F_0:\ldots: F_k]$ as a meromorphic map in $\P^k.$
 For $d\geq 2,$ let $\Fc_d(\P^k)$ be the  space of singular holomorphic  foliations of degree $d$ in $\P^k.$
 Using  the above  form of $Z,$  we can show that  $\Fc_d(\P^k)$  can be canonically identified  with a Zariski open subset of $\P^N,$ where
 $N:=(d+k+1){(d+k-1)!\over (k-1)! d!} -1 $ (see \cite{Brunella}).  The  next result describes  the typical properties of a generic foliation $\Fc\in\Fc_d(\P^k).$
 
 \begin{theorem}\label{T:generic} Let $d,k>1$.

  \begin{enumerate}
  \item {\rm (Jouanolou \cite{Jouanolou}, Lins Neto-Soares \cite{NetoSoares}).}
  There is a real Zariski dense open set $\mathcal H(d)\subset  \Fc_d(\P^k)$ such that for every $\Fc\in \mathcal H(d),$
  all  the singularities of $\Fc$ are    hyperbolic and  $\Fc$ do not possess any invariant algebraic  curve.
  \item       {\rm (Glutsyuk  \cite{Glutsyuk},  Lins Neto \cite{Neto}).}
  If all  the singularities of a  foliation $\Fc\in \Fc_d(\P^k)$ are    non-degenerate, then $\Fc$ is hyperbolic.
   \item  {\rm (Brunella \cite{Brunella}).}   If 
  all  the singularities of a  foliation $\Fc\in \Fc_d(\P^k)$ are    hyperbolic and  $\Fc$ do not possess any invariant algebraic  curve, then $\Fc$ admits no  nontrivial  directed positive closed current.
  \end{enumerate}
\end{theorem}
  
  Moreover, Loray-Rebelo \cite{LorayRebelo}  constructed a nonempty  open  set  $\mathcal U(d)$ of $\Fc_d(\P^k)$
   such that every leaf of $\Fc\in \mathcal U(d)$ is  dense.
   
   Now  we come to the notion of directed positive harmonic currents on singular Riemann surface  laminations.
   \begin{definition}\label{D:Directed_hamonic_currents_with_sing}  {\rm  (Berndtsson-Sibony \cite{BerndtssonSibony}, Forn{\ae}ss-Sibony \cite{FornaessSibony05,FornaessSibony08}).}  \rm 
   Let  $(X,\Lc,E)$   be  a  Riemann surface lamination  with  singularities, where 
$X$ is a closed
subset of a
complex manifold $M$  and    the leaves of $(X\setminus E,\Lc)$ are Riemann surfaces holomorphically immersed in
$M$. 
A {\it directed  harmonic  current} on $(X,\Lc,E)$ is  a  positive harmonic  current $T$ of bidimension $(1,1)$ on $M$ such that  the  support of $T$  is contained in $X$ and  that
the restriction of $T$  on $X\setminus E$ is  a  directed harmonic  current on the Riemann surface lamination $(X,\Lc)$ in the sense of Definition \ref{D:Directed_hamonic_currents}.
   \end{definition}
   The  existence of directed  positive harmonic  currents for  compact (nonsingular) laminations was   proved by Garnett \cite{Garnett}.
   The case  of compact singular Riemann surface laminations was  proved by 
   Berndtsson-Sibony under  reasonable assumptions.
   \begin{theorem}\label{T:existence_harmonic_currents} {\rm  (Berndtsson-Sibony \cite{BerndtssonSibony}, see also  \cite[Theorem 23]
   {FornaessSibony08}).}
   Let $(X,\Lc,E)$ be a singular  Riemann surface lamination as in the assumption  of Definition  \ref{D:Directed_hamonic_currents_with_sing}.
   Assume moreover that $X$ is  compact  and  $E$ is locally pluripolar in $M.$ Then there is a nonzero directed positive harmonic  current $T.$ 
   In particular,  if the set $E$ does not support any nonzero positive harmonic current (e.g. if  $\Lambda_2(E)=0,$ where $\Lambda_2$ denotes the two dimensional Hausdorff measure),
   then the  restriction of such   a current $T$
   on $X\setminus E$ induces 
     a nonzero directed positive harmonic  current on $(X\setminus E,\Lc).$  
   \end{theorem}
   When  a leaf $L_x$  is   hyperbolic, an average on $L_x$   was introduced by Forn{\ae}ss-Sibony \cite{FornaessSibony05} (see also  \cite[Corolary 3]{FornaessSibony08}). 
   It allows another construction of  directed positive  harmonic currents. 
 By  Theorem \ref{T:generic}, Theorem \ref{T:existence_harmonic_currents} applies to every generic foliation in  $\P^k$ with a given degree $d>1.$  
 
 \subsection{General laminations/foliations}\label{SS:general_notions}
We  formulate  some  general notions of  laminations and foliations. Although  they are not  the main  topic  of this  article,
some  results  presented  here could  be   extended to  these general objects.  Let $l\geq1$ be an integer.
 
 An {\it   $l$-dimensional  lamination}     $(X,\Lc)$   is  the  data of a locally compact space $X$ and  a {\it (lamination)  atlas} $\Lc$ 
of   with (laminated) charts 
$$\Phi_p:\U_p\rightarrow \B_p\times \T_p.$$
Here, $\T_p$ is a locally compact metric space, $\B_p$ is a domain in $\R^l$,  $\U_p$ is  an open set in 
$X,$ and  
$\Phi_p$ is  a homeomorphism,  and  all the changes of coordinates $\Phi_p\circ\Phi_q^{-1}$ are of the form
$$x=(y,t)\mapsto x'=(y',t'), \quad y'=\Psi(y,t),\quad t'=\Lambda(t),$$
 where $\Psi,$ $\Lambda$ are continuous functions. 
 We call {\it  $l$-dimensional  lamination with singularities} the data
$(X,\Lc,E)$ where $X$ is  a locally compact space, $E$ is a closed
subset of $X$ and $(X\setminus E,\Lc)$ is an $l$-dimensional    lamination.
 The set $E$ is {\it the singularity set} of the lamination.
 
 An {\it   $l$-dimensional  Riemannian  foliation}     $(X,\Lc)$   is  the  data of  a {\it (lamination)  atlas} $\Lc$ 
of  a Riemannian  manifold $X$ of dimension $k\geq l$ with (laminated) charts 
$$\Phi_p:\U_p\rightarrow \B_p\times \T_p.$$
Here, $\T_p$ is a domain in $\R^{k-l},$ $\B_p$ is a domain in $\R^l$,  $\U_p$ is  an open set in 
$X,$ and  
$\Phi_p$ is  a homeomorphism,  and  all the changes of coordinates $\Phi_p\circ\Phi_q^{-1}$ are of the form
$$x=(y,t)\mapsto x'=(y',t'), \quad y'=\Psi(y,t),\quad t'=\Lambda(t),$$
 where $\Psi,$ $\Lambda$ are smooth functions. 
We call {\it  $l$-dimensional  Riemannian  foliation with singularities} the data
$(X,\Lc,E)$ where $X$ is  a Riemannian  manifold $X$ of dimension $k\geq l$, $E$ is a closed
subset of $X$ and $(X\setminus E,\Lc)$ is an $l$-dimensional  Riemannian  foliation.
 The set $E$ is {\it the singularity set} of the foliation.

 By considering classical constructions such as projective limits and suspensions, 
Forn{\ae}ss-Sibony-Wold  obtain in \cite{FornaessSibonyWold}  examples of laminations possessing many, or few, directed positive closed or directed positive harmonic currents. 
 They  also give examples of compact  laminations by complex manifolds of dimension $\geq 2$  with no nonzero directed positive  harmonic current. 
 This is in contrast with the Riemann surface case, where such a current always exists.
 
\section{Regularity of the leafwise Poincar\'e metric} \label{S:Poincare}

Let $(X,\Lc,E)$ be a  hyperbolic Riemann surface  lamination  with singularities. Let $g_P$ be    the leafwise Poincar\'e metric for the  lamination  $(X\setminus  E,\Lc)$ given in Subsection  \ref{SS:Poincare}.
Let $g_X$ be a Hermitian metric on the leaves which is transversally smooth. We can construct such a metric on flow boxes and glue them using a partition of unity.
We have
\begin{equation}\label{e:eta} g_X=\eta^2g_P\quad \mbox{where}\quad 
\eta(x):=\|D\phi_x(0)\|.
\end{equation}
Here, $\phi_x$ is  defined in   \eqref{e:covering_map}, and for the norm of the differential $D\phi_x$ we use the Poincar\'e metric on $\D$ and the Hermitian metric $g_X$ on $L_x$. 

The extremal property of the Poincar\'e metric implies that
$$\eta(x)=\sup\big\{\|(D\phi)(0)\|,\quad \phi:\D\to L \mbox{ holomorphic such that } \phi(0)=x\big\}.$$
Using a map sending $\D$ to a plaque, we see that the function $\eta$ is locally bounded from below on $X\setminus E$ by a strictly positive constant. 
When $X$ is compact and $E=\varnothing,$ the classical Brody lemma (see \cite[p.100]{Kobayashi}) implies that $\eta$ is also bounded from above.

The continuity of the function  $\eta$ was studied by Candel, Ghys, Verjovsky,   see \cite{Candel,Ghys,Verjovsky}. The survey \cite{FornaessSibony08} establishes the result as a consequence of Royden's lemma. 
Indeed with his lemma, Royden proved the upper-semicontinuity of the infinitesimal Kobayashi metric in a Kobayashi hyperbolic manifold (see \cite[p.91 and p.153]{Kobayashi}).
The  following theorem gives refinements of  the previous results.

\begin{theorem} \label{T:Poincare}{\rm (Dinh-Nguyen-Sibony \cite{DinhNguyenSibony14a}). }
Let $(X,\Lc)$ be a transversally smooth compact lamination by hyperbolic Riemann surfaces.
Then the Poincar\'e metric on the leaves  is H\"older continuous, that is, the function $\eta$ defined in \eqref{e:eta} is H\"older continuous on $X$. 
Moreover, the exponent of   H\"{o}lder continuity can be estimated in geometric terms.
\end{theorem}

The main tool of the proof of Theorem  \ref{T:Poincare}  is
to use  Beltrami's equation in order to
compare universal covering maps of  any leaf $L_y$ near a given leaf $L_x$.
More precisely, 
for  $R>0$ let $\D_R$ be  the  disc of center $0$ with  radius $R$  with respect to  the  Poincar\'e metric on $\D$ (see Main Notations in Section \ref{S:introduction}). 
We first  construct a non-holomorphic parametrization $\psi$ from $\D_R$ to $L_y$ which is close to a universal covering map 
$\phi_x:\D\to L_x$ for all $R$ large enough. Next, precise geometric estimates on $\psi$ allow us to modify it, using Beltrami's equation. 
We then obtain a holomorphic map that we can explicitly compare with a universal covering map $\phi_y:\D\to L_y$.

 Next,    we investigate the regularity  of the leafwise Poincar\'e metric $g_P$  of a compact singular holomorphic  foliation. Here  
  an important difficulty   emerges: a leaf of the foliation may visit singular flow boxes without any obvious rule. 
  We are interested  in the following class of laminations.
 
\begin{definition}\label{D:uniform_hyperbolic_laminations} {\rm (Dinh-Nguyen-Sibony \cite{DinhNguyenSibony14b}). }  \rm
 A  hyperbolic Riemann surface  lamination  with singularities   $(X,\Lc,E)$  with $X$ compact
 is  said  to be {\it Brody hyperbolic}  if  there is  a constant  $c_0>0$  such that
$$\|D\phi(0)\| \leq  c_0$$ 
for all  holomorphic maps  $\phi$ from $\D$  into  a leaf.
\end{definition}

 \begin{remark}\rm \label{R:Brody_sufficiency} 
It is clear that if the lamination is Brody hyperbolic then its leaves are hyperbolic in the sense of Kobayashi.
Conversely, the  Brody hyperbolicity is  a consequence of
the non-existence  of  holomorphic non-constant maps
$\C\rightarrow X$  such that out of  $E$ the image of $\C$ is  locally contained in leaves, 
see \cite[Theorem 15]{FornaessSibony08}.

On the other hand,   Lins Neto proved  in   \cite{Neto} that  for every   holomorphic foliation of degree larger than 1  in $\P^k$, with non-degenerate singularities,  
there is a smooth metric with negative curvature on its tangent bundle, see also Glutsyuk \cite{Glutsyuk}.
Hence, these foliations are Brody hyperbolic.
Consequently,  holomorphic foliations  in $\P^k$  are generically  Brody hyperbolic, see Theorem \ref{T:generic} (1).
\end{remark}

 Denote by   $\lof(\cdot):=1+|\log(\cdot)|$ a log-type function, and  by $\dist$  the distance on $X$ induced by  the Hermitian metric $g_X.$
The following result is   a  counterpart of Theorem \ref{T:Poincare} in the context of singular holomorphic foliations. 

\begin{theorem} \label{T:Poincare_bis} {\rm (Dinh-Nguyen-Sibony \cite{DinhNguyenSibony14b}). } 
Let $(X,\Lc,E)$ be a Brody hyperbolic singular holomorphic foliation  on a Hermitian compact complex manifold $X$. Assume that the singular set $E$  is 
finite  and  that all  points of $E$  are  linearizable. 
 Then, there are constants  $c>0$ and $0<\alpha<1$   such that
 $$|\eta(x)-\eta(y)|\leq  c\Big ( {\max\{\lof \dist(x,E), \lof\dist(y,E)\}\over \lof\dist(x,y)}\Big)^\alpha$$
 for all $x,y$ in $X\setminus E$. 
\end{theorem}

  To prove  this  theorem, we analyze the behavior and get an explicit estimate on the modulus of continuity of the Poincar\'e metric on leaves.
The following estimates are crucial in our method.  They  are also  useful in other  problems.

\begin{proposition} \label{P:Poincare}{\rm (Dinh-Nguyen-Sibony \cite{DinhNguyenSibony14b}). }
Under the hypotheses
of Theorem  \ref{T:Poincare_bis},  there exists  a  constant $c_1>1$  such that 
$$c_1^{-1} s \lof s \leq\eta(x)  \leq c_1  s \lof s$$
for $x\in X\setminus E$  and $s:=\dist(x,E)$.  
\end{proposition}

The  Poincar\'e metric on  the leaves of a hyperbolic  foliation is a fundamental  object  which is
extremely delicate  to understand.   As  we  see in Theorem \ref{T:Poincare_bis}, the  regularity in the direction  transverse to the foliation is  quite  weak.
This is  partly  due to   the presence  of the  singularities.
We  end the section with the  following open  question.
\begin{problem}  \rm  
Let $(X,\Lc,E)$ be a compact  singular  holomorphic  foliation by hyperbolic Riemann surfaces.
Assume that every point $a\in E$ is   a non-degenerate singularity. 
Study the regularity  of the function $\eta.$
In case  $\dim(X)=2,$  we may investigate the  problem   where the singularities are not necessarily non-degenerate. 
\end{problem}

 \section{Mass-distribution of directed positive harmonic currents}\label{S:Mass}

Let $(X,\Lc,E)$ be a  singular holomorphic  foliation and let
 $T$ be a positive harmonic current on $X\setminus E$. By  Proposition  \ref{P:extension} its mass with
respect to any Hermitian metric on $X$ is  finite. 
We call {\it Poincar{\'e}  mass} of $T$ the mass of $T$ with respect to
Poincar{\'e}  metric $g_P$ on $X\setminus E$, i.e. the mass of the
positive measure $m_P:=T\wedge g_P$.
A priori, Poincar{\'e}  mass may be
infinite near the singular points. The following proposition
gives us a criterion for the finiteness of this mass. It can be applied
to generic foliations in $\P^k$  (see Theorem  \ref{T:generic}).
 
\begin{proposition} \label{P:Poincare_mass} {\rm (Dinh-Nguyen-Sibony \cite{DinhNguyenSibony12}).}
Let $(X,\Lc,E)$ be a singular  holomorphic  foliation. If $a\in E$ is a linearizable singularity, 
then any positive harmonic current on $X$ has locally finite
Poincar{\'e}  mass near $a$.
\end{proposition}

The proof of this  result is based on the finiteness of the Lelong number of $T$  at $a$  (see Proposition \ref{P:Skoda}).
In dimension $2$ we have  a more precise result when   $T$ is  directed and the singular point is  hyperbolic. 
\begin{theorem} \label{T:Lelong}
{\rm (Nguyen \cite{NguyenVietAnh17c}).} Let $(X,\Lc,E)$ be a singular  holomorphic  foliation   with $\dim X=2.$ If $a\in E$ is a hyperbolic singularity, 
then for any directed positive harmonic current $T$ on $X$  which does  not give mass to any of the two separatrices at $a,$    the   Lelong  number of $T$ at $a$  vanishes.  
 \end{theorem}
An  immediate consequence  of Theorem \ref{T:Lelong} is the following result on the Lelong numbers of  a directed harmonic current.  
\begin{corollary}\label{C:Lelong}
Let $\Fc=(X,\Fc,E)$    be  
a     singular holomorphic foliation   with $X$ a    compact complex surface. 
Assume that all the singularities  are hyperbolic and  that the foliation  has no invariant analytic  curve.
Then for every   harmonic  current $T$  directed by $\Fc,$   the   Lelong  number of $T$ vanishes everywhere   in $X.$  
 \end{corollary}

The  above   corollary can be   applied  to  every generic  foliation in $\P^2$ with a given degree $d>1$
(see Theorem  \ref{T:generic}).

We can apply  Corollary \ref{C:Lelong} to   study the   recurrence  of a  generic leaf. More specifically,
let $T$ be a positive harmonic current directed by a singular holomorphic  foliation $(X,\Lc,E)$ with $X$  a compact  complex  surface.   
Assume that all the singularities  are hyperbolic and  that the foliation  has no invariant analytic  curve.
Consider  the positive measure $m_P:=T\wedge g_P.$ We know  by Proposition  \ref{P:Poincare_mass} that $m_P$ is a  finite measure.
Given    a point $x\in X$ and  a  $m_P$-generic  point $a\in X\setminus  E,$ we want to know
how  often   the  leaf  $L_a$ visits  the ball $B(x,r)$   as  $r\searrow 0.$
Here  $B(x,r)$ 
denotes the open 
ball with  center $x$ and radius $r$ with respect to a fixed metric on $X.$

 Let us  introduce  some more  notation and terminology.
Denote by
$r\D$ the disc of center $0$ and of radius $r$ with $0<r<1$. In the
Poincar{\'e} disc $(\D,\omega_P),$ 
$r\D$ is also the disc of center 0 and of radius 
\begin{equation}\label{e:radii_conversion}
R:=\log{1+r\over 1-r}\cdot
\end{equation}
So, we will also denote by $\D_R$ this disc, and by $\partial \D_R$ its  boundary  which is also the Poincar\'e circle of center $0$ and radius $R.$

Together with Dinh and Sibony,  we introduce  in  \cite{NguyenVietAnh17c} the following  indicator.  
\begin{definition}\label{D:visibility}\rm For  each $r>0,$ 
the {\it visibility  of  a point $a\in X\setminus E$  within   distance $r$  from a  point $x\in M$} is the number
$$N(a,x,r)= \limsup_{R\to \infty}\frac{1}{R}\int\limits_0^R \Big( \int_{\theta=0}^1  \textbf{1}_{ B(x,r)}\big (\phi_a(s_te^{2\pi i\theta})\big)d\theta\Big) dt 
 \in [0,1],$$
where $ \textbf{1}_{ B(x,r)}$ is
the  characteristic  function  associated to the set $B(x,r),$ and $s_t$  is  defined by  the  relation $t=\log{1+s_t\over 1-s_t}\cdot,$ that is,  $s_t\D=\D_t .$
\end{definition}
Geometrically, $N(a,x,r)$  is the  average, as $R\to\infty,$  over  the hyperbolic time $t\in [0, R]$ of  the Lebesgue measure of the set
$\{  \theta\in[0,1]:\  \phi_a(s_te^{2\pi i\theta})  \in B(x,r)\}.$ 
The last  quantity  may be interpreted  as the portion  which hits  $B(x,r)$ of the Poincar\'e circle of radius $t$  with center $a$ spanned on the leaf $L_a.$  


We combine  Corollary  \ref{C:Lelong} and the  so-called  geometric  ergodic  theorem  (Theorem \ref{T:geometric_ergodic}) which will be presented in the  next section.
Consequently, we  obtain the  following  upper bound on the visibility of a  generic point. 
\begin{theorem}\label{T:visibility}
 We keep  the  above  hypothesis and notation.  Then  for $m_P$-almost  every  point $a\in X\setminus E$
 and for every point $x\in X,$  we have that
 $$
 N(a,x,r)=\begin{cases}
 o( r^2), & x\in X\setminus E;\\
 o(|\log r|^{-1}), & x\in E.
 \end{cases}
 $$
\end{theorem}

Here  are some open questions.

\begin{problem}\rm
Let $\Fc=(X,\Lc,E)$ be a compact singular  holomorphic  foliation by hyperbolic Riemann surfaces. Let  $a\in E$ be a non-degenerate singularity.
Find sufficient conditions on $\Fc$ and $a$   so that   
 any positive harmonic current on $X$ has locally finite
Poincar{\'e}  mass near $a$.
\end{problem}

  \begin{problem}\rm
    Can one  generalize      Theorem \ref{T:Lelong}    to higher dimensions? 
  \end{problem}

  \begin{problem}\rm
   Find  an effective  lower bound   on the  visibility  of a generic point.  
  \end{problem}

\section{Heat  equation and  ergodic theorems} \label{S:Ergodic_theorems}

In this  section  we  will report some  recent  techniques used   to obtain  ergodic  theorems for Riemann surface laminations.
\subsection{Ergodic theorems  associated with the  heat  diffusions}\label{SS:heat_diffusions}

In collaboration with  Dinh and Sibony  \cite{DinhNguyenSibony12},  we introduce the heat equation relative to a positive harmonic closed current
and apply it to the directed positive harmonic currents of a  Riemann surface
laminations  with singularities.  This permits to construct the heat diffusion with respect to
various Laplacians that could be defined
almost everywhere with respect to the positive harmonic current.

More concretely, let  $(X,\Lc,E)$   be  a compact Riemann surface lamination  with  singularities, where 
$X$ is a (compact)
subset of a
complex manifold $M$  and    the leaves of $(X\setminus E,\Lc)$ are Riemann surfaces holomorphically immersed in
$M$. For simplicity, fix a Hermitian form $g_M$ on $M$. Let $T$ be a directed positive  harmonic  current. So,
$T\wedge g_M$ is a positive measure. Consider a positive $(1,1)$-form $\beta$ which is defined
almost everywhere on $M$ with respect to $T\wedge g_P.$ We say that $T$ is {\it $\beta$-regular} if 
$T\wedge\beta$ is of finite mass and  $i\tau\wedge \overline{\tau}\wedge T\leq  T\wedge\beta,$ 
where $\tau$
is a $(1,0)$-form 
defined almost everywhere with respect to $T\wedge g_P$ such that
$\partial T=\tau\wedge T .$ Under  the notation of  Proposition \ref{P:decomposition},  we  see that  $\tau= h_t^{-1}\partial h_t$ on the  plaque passing  through $t\in\T$
for $\nu$-almost every $t\in\T$. The following result  gives a typical example of $\beta$-regularity with $\beta:=g_P.$

\begin{proposition}\label{P:g_P-regular}
Let $(X,\Lc,E)$ be a lamination  as  above.
Let $\Fc$ be  a  singular holomorphic  foliation on $M$ such that the restriction of $\Fc$ on $X\setminus E$  induces $\Lc$ and that
all points of $E$ are linearizable  singularities of $\Fc.$  
Then every directed harmonic current $T$ on $(X,\Lc,E)$  is $g_P$-regular. 
\end{proposition}
 Using Proposition \ref{P:Poincare_mass},  Proposition \ref{P:g_P-regular}  follows  essentially  from  Proposition 3 and  Theorem 9 in  \cite{FornaessSibony10}.

 We define  the Laplacian $\Delta_\beta$ by
\begin{equation}\label{e:Laplacian}
(\Delta_\beta u)T \wedge \beta :=\ddc u \wedge T\mbox{  for  } u\in\Cc^\infty_0(M).
 \end{equation}
We  will extend  the  definition of $\Delta_\beta$  to larger  spaces, suitable  for developing   $L^2$-techniques. To this end, 
let $m_\beta$ denote the measure $T\wedge\beta$ and consider  the  Hilbert space $L:=L^2(m_\beta).$ 
We also  introduce  the Hilbert space  $H=H^1_\beta(T)\subset L^2(m_\beta)$ 
associated with $T$ and  $\beta$  as the  completion of $\Cc^\infty_0(M)$ with  respect  
to the Dirichlet norm
$$\|u\|^2_{H^1_\beta}:=\int |u|^2  T\wedge\beta +i\int \partial u
\wedge \dbar u\wedge T.$$
 Using  the  assumption that $T$ is $\beta$-regular, we can show that there exists a semi-group of contractions $S(t) : L \to L,$
$t \in \R^+$ such that for every  function $u_0\in H,$   $u(t, \cdot) := S(t)u_0$ satisfies
$${\partial u(t, \cdot)\over\partial t}
=\Delta_\beta u(t, \cdot) \mbox{  and }  u(0, \cdot) = u_0.$$
Recall that
a family $S(t):L\rightarrow L$, $t\in\R_+$, is {\it a semi-group of
  contractions} if $S(t+t')=S(t)\circ S(t')$ and if $\|S(t)\|\leq 1$
for all $t,t'\geq 0$.

To prove  this  result  we  use  functional  analysis  (Hille-Yosida theorem, Lax-Milgram theorem  etc.).
We also   use  Stokes' theorem on $M.$
It is  worthy noting that  Garnett \cite{Garnett} and   Candel \cite{Candel2} also   solve the  heat equation. But they consider the  case without singularities. Moreover,
they  solve the equation  pointwise, that is, in the space of smooth  functions. So their methods  are  quite diffrent from  
ours. Indeed,  we solve  the  equation with respect to a  harmonic  current,  in a suitable $L^2$-space. 

The  following result is   an ergodic theorem associated to the heat
diffusions.
\begin{theorem}\label{T:diffusions} {\rm (Dinh-Nguyen-Sibony \cite{DinhNguyenSibony12}). }
 We keep the above hypothesis and notation.
  Then 
\begin{enumerate}
 \item 
the  
measure $m_\beta$ is  $S(t)$-invariant (that is,  $\langle S(t)u,m_\beta\rangle=\langle u,m_\beta\rangle $ for every $u\in L$), and 
$S(t)$ is a positive contraction in $L^p(m_\beta)$ for all $1\leq p\leq\infty$ (that is, $\|S(t)u\|_{L^p(m_\beta)}\leq  \|u\|_{L^p(m_\beta)}$ for every $u\in L$); 
\item
for all $u_0\in  L^p(m_\beta),$  $1\leq p<\infty$, the average 
$$\frac{1}{R}\int_0^R S(t)u_0 dt$$ 
converges pointwise  $m_\beta$-almost everywhere
 and also in $L^p(m)$ to 
an $S(t)$-invariant function  $u_0^*$ when $R$ goes to infinity.
Moreover, $u_0^*$ is constant on the leaf $L_a$ for $m_\beta$-almost every $a$.
If $m_\beta$ is an extremal harmonic measure, then $u$ is constant $m_\beta$-almost everywhere.
\end{enumerate}
\end{theorem}

Combining  Proposition \ref{P:g_P-regular} and Theorem \ref{T:diffusions}, we obtain the following 
relation between   harmonic measures and directed positive harmonic currents which is a complement to  Theorem \ref{thm_harmonic_currents_vs_measures} (ii).
\begin{proposition}\label{P:harmonic_currents_vs_measures}
Let $(X,\Lc,E)$ be a lamination  as  above.
Let $\Fc$ be  a  singular holomorphic  foliation on $M$ such that the restriction of $\Fc$ on $X\setminus E$  induces $\Lc$ and that
all points of $E$ are linearizable  singularities of $\Fc.$  
Then the map $T\mapsto  \mu=T\wedge g_P$   is  a bijection   from the  convex  cone of  all  directed positive harmonic  currents $T$  onto   
the convex  cone of all harmonic  measures $\mu$.
\end{proposition}
\subsection{Geometric ergodic theorems}

In this subsection, we will give an analogue  of  
Birkhoff's ergodic theorem in the  context of a compact  Riemann surface  lamination $(X,\Lc,E)$
with  singularities. 
 Our  ergodic theorem is of  geometric nature  and   it is close to Birkhoff's averaging on orbits of a
dynamical system.   Here  the averaging is  on  hyperbolic  leaves and the  time is the  hyperbolic time.  

Let $(X,\Lc,E)$ be a  Riemann surface  lamination with singularities which is  embedded in a complex manifold $M$ as in Subsection \ref{SS:heat_diffusions}.
Let $T$ be a directed positive
harmonic current on $(X,\Lc,E)$ such that $T$ is $g_P$-regular.
A leaf $L_x$ is  called {\it  parabolic} if it is  not  hyperbolic. We
assume that $T$ has no mass on the union of parabolic leaves and that
$m_P:=T\wedge g_P$ is a probability measure. So by Proposition \ref{thm_harmonic_currents_vs_measures}, $m_P$ is a quasi-harmonic measure on $X$ with respect to $g_P$.

For any point $x\in X\setminus E$ such that the corresponding leaf
$L_x$ is hyperbolic,   
let $\phi_x:\D\rightarrow L_x$  be given by  \eqref{e:covering_map}. Denote by
$r\D$ the disc of center 0 and of radius $r$ with $0<r<1$.  Recall  from \eqref{e:radii_conversion} that in the
Poincar{\'e} metric, this is also the disc of center 0 and of radius 
$ R:=\log{1+r\over 1-r},$
and we will also denote by $\D_R$ this disc.
For all $0<R<\infty$, consider 
\begin{equation}\label{e:m_x,R}\begin{split}  m_{x,R}&:=\frac{1}{M_R}(\phi_x)_* \big(\log^+ \frac{r}{|\zeta|}g_P\big),\\
 \tau_{x,R}&:=\frac{1}{M_R}(\phi_x)_* \big(\log^+ \frac{r}{|\zeta|}\big).
 \end{split}
 \end{equation}
where $\log^+:=\max\{\log,0\}$, $g_P$ denotes also the Poincar{\'e} metric on $\D$ and 
\begin{equation*}
M_R:= \int \log^+ \frac{r}{|\zeta|} g_P=\int \log^+ \frac{r}{|\zeta|} \frac{2}{(1-|\zeta|^2)^2}
id\zeta\wedge d\overline\zeta.
\end{equation*}
So, $m_{x,R}$ (resp. $\tau_{x,R}$) is a probability measure  (resp.  a directed  positive  current of bidimension $(1,1)$) which depends on $x,R;$ but 
does not depend on the choice of $\phi_x$.

\begin{theorem} \label{T:geometric_ergodic}{\rm (Dinh-Nguyen-Sibony \cite{DinhNguyenSibony12}). }
 We  keep the above hypothesis and notation.
Assume in addition that the current  $T$ is extremal.    
Then for almost every point $x\in X$ with respect to the measure
$m_P:=T\wedge g_P$, the measure $m_{x,R}$ defined above converges to $m_P$
when $R\to\infty$.  
\end{theorem}

 To  prove  the theorem,  our main ingredient  is a  delicate estimate on the heat kernel of the Poincar\'e disc
 (see   \cite[p. 370, line 8-]{DinhNguyenSibony12} for its statement).
 This estimate allows us to deduce  the the  desired result  from 
  the ergodic theorem associated to the heat
diffusions (Theorem \ref{T:diffusions}).

 \begin{remark}
  \rm
Let $(X,\Lc)$ be a compact lamination by Riemann surfaces without singularities.
Let $T$ be a positive harmonic current directed by the lamination which is
extremal, with full mass on hyperbolic leaves and with Poincar\'e mass 1. Then, the conclusions of Theorem \ref{T:diffusions} and
Theorem \ref{T:geometric_ergodic} are still valid. The proofs are
essentially the same but we need  to use  a finite partition of unity  for $X$  instead of  applying  Stokes' theorem for $M.$
Moreover, Theorem  \ref{T:diffusions} still holds for compact smooth $p$-dimensional  laminations in the sense of Subsection  \ref{SS:general_notions}.
 \end{remark}

\subsection{Unique ergodicity theorems}

In  \cite{FornaessSibony05}   Forn{\ae}ss  and Sibony develop the theory of  harmonic  currents of finite energy.
They introduce  a notion of energy for positive harmonic  currents of bidegree $(1,1)$  on a compact K\"ahler manifold $(X,\omega)$
of dimension $k\geq 2.$
This  allows to define  $\int_X T\wedge T\wedge \omega^{k-2}$ for  every positive harmonic  current $T$ of bidegree $(1,1)$ on $X.$
This theory applies to directed  positive harmonic  currents on singular holomorphic  foliations on compact K\"ahler surfaces.

In  \cite{FornaessSibony05, FornaessSibony10}   Forn{\ae}ss  and Sibony also  develop a geometric intersection theory 
for directed  positive harmonic  currents on singular holomorphic  foliations on $\P^2.$

Combining these two theories, they obtain  the following remarkable unique   ergodicity result for   singular holomorphic foliations
without  invariant algebraic  curves.

\begin{theorem} \label{T:FornaessSibony}{\rm (Forn{\ae}ss-Sibony \cite{FornaessSibony10}). } Let $\Fc$ be a singular  holomorphic foliation  in $\P^2$ 
whose singularities are all hyperbolic. Assume that $\Fc$ has no invariant
algebraic curve. Then $\Fc$  has a unique directed positive  harmonic  
current of mass $1.$ Moreover, this unique current $T$ is not closed.
In particular,  for every point $x$ outside the singularity set of $\Fc,$
the current $\tau_{x,R}$ defined in \eqref{e:m_x,R} converges to $T$
when $R\to\infty$.  
\end{theorem}

The  case  where $\Fc$  possesses  invariant algebraic  curves has recently  been answered.

\begin{theorem} \label{T:Dinh-Sibony} {\rm (Dinh-Sibony \cite{DinhSibony17a}). }
Let $\Fc$ be a singular  holomorphic foliation  in $\P^2$ 
whose singularities are all hyperbolic. Assume that
 $\Fc$
admits a finite number of invariant algebraic curves. Then   any directed   positive harmonic current is a linear combination of the currents of integration on
these curves.  In particular, all directed   positive harmonic currents
are closed.
\end{theorem}

 Theorem \ref{T:Dinh-Sibony} is  suprising even in the special  case where $\Fc$  admits the line at infinity $L_\infty$  as an invariant curve.
 Let $\Fc$ be a generic  foliation of a given degree $d>1$ with this property.
 By Khudai-Veronov 
\cite{IlyashenkoYakovenko}, all  leaves  (except $L_\infty$) of $\Fc$ are dense. So by intuition from  Theorem \ref{T:FornaessSibony} one could  expect   that there  should be 
a directed harmonic  current with support $\P^2.$ However, Theorem \ref{T:Dinh-Sibony}  says that this  intuition is   false.

To prove Theorem   \ref{T:Dinh-Sibony} we need to to show that if $T$ is a positive harmonic
current directed by $\Fc$ having no mass on any leaf, then $T$ is zero.
For this purpose,   Dinh and Sibony \cite{DinhSibony17a} develop a theory of densities of positive $\ddc$-closed $(1,1)$-currents in a compact K\"ahler surface. A related
theory was previously developed by these authors in \cite{DinhSibony12} for positive closed currents  defined on  compact K\"ahler manifolds. 
Applications of these  theories in   complex dynamics of higher dimension could be found in \cite{DinhNguyenTruong,DinhSibony16, DinhSibony17b} etc.

  \begin{problem}\rm
    Are there  any versions of     Theorem \ref{T:FornaessSibony} and    Theorem   \ref{T:Dinh-Sibony}  for
 compact K\"ahler surfaces  ? 
  \end{problem}

\section{Topological and  metric entropies for hyperbolic Riemann surface laminations} \label{S:Entropy}


 Ghys-Langevin-Walczak  introduced  in \cite{GhysLangevinWalczak}   a notion  of geometric entropy  for compact transversally smooth Riemannian foliations $(X,\Lc)$
 (see also Candel-Conlon  \cite{CandelConlon1}  and Walczak  \cite{Walczak} for  recent expositions). They  prove that this entropy is always finite. In fact,  their notion is
related to the entropy of the holonomy pseudogroup, which depends on the  chosen  generators. It also  depends on the choice of the metric on the  ambient manifold $X.$
The basic idea
is to quantify how  much leaves get far apart transversally. The transverse regularity of the metric on leaves and the lack of singularities  play a role in the  finiteness of the entropy.

In  \cite{DinhNguyenSibony14a} Dinh-Nguyen-Sibony  introduce a general notion of entropy, which permits to
 describe some natural situations in dynamics and in laminations/foliation theory. 
This new notion of entropy   contains a large number of classical situations.
In particular, it also  includes  Riemannian foliations  with singularities. 
Another interesting fact is that this entropy is related to an increasing family of distances as in  Bowen's point of view \cite{Bowen}.
This allows,  for example,  for introducing other dynamical notions like metric entropy, local entropies etc. 


Let $X$ be a metric space endowed with a distance $\dist_X$. Consider a family $\Dc=\{\dist_t\}$ of distances on $X$ indexed by $t\in \R^+$. We can 
also replace $\R^+$ by $\N$ and in practice we often have that $\dist_0=\dist_X$ and that 
$\dist_t$ is  increasing  with respect to  $t\geq 0$.

Let $Y$ be a non-empty subset of $X$. Denote by $N(Y,t,\epsilon)$ the minimal number of balls of radius $\epsilon$ with respect to the distance $\dist_t$ needed to cover $Y.$
Define the {\it  entropy} of $Y$ with respect to $\Dc$ by
\begin{equation}\label{e:entropy}
h_\Dc(Y):=\sup_{\epsilon>0} \limsup_{t\rightarrow\infty} {1\over t} \log  N(Y,t,\epsilon).
 \end{equation}
When $Y=X$ we will denote by $h_\Dc$ this entropy. 

Observe that when $\dist_t$ is increasing,  $N(Y,t,\epsilon)$ is  increasing  with respect to  $t\geq 0$. Moreover,   
$$ \limsup_{t\rightarrow\infty} {1\over t} \log  N(Y,t,\epsilon)$$  is  increasing
when $\epsilon$  decreases. So, in the above definition, we can replace $\sup_{\epsilon>0}$ by $\lim_{\epsilon\to 0^+}$. 
If $\Dc=\{\dist_t\}$  and $\Dc'=\{\dist'_t\}$ are two families of distances on $X$  such that  $\dist'_t\geq  A\dist_t$ for all $t$ with  a fixed constant $A>0$,  then  $h_{\Dc'}\geq h_\Dc$.

A subset $F\subset X$ is said to be {\it $(t,\epsilon)$-separated} if for all distinct points $x,y$ in $F$ we have $\dist_t(x,y)>\epsilon$.
Let $M(Y,t,\epsilon)$ denote the maximal number of points in a $(t,\epsilon)$-separated family $F\subset Y$. 
We   record  here a  simple  relation between $N(Y,t,\epsilon)$ and $M(Y,t,\epsilon).$

\begin{lemma}\label{L:comparison} 
We have
$$N(Y,t,\epsilon)\leq   M(Y,t,\epsilon)\leq N(Y,t,\epsilon/2).$$
\end{lemma}
An important   consequence of Lemma \ref{L:comparison} is that we can formulate  the entropy of a subset $Y\subset X$  using
$M(Y,t,\epsilon)$ instead of $N(Y,t,\epsilon):$
$$h_\Dc(Y)= \sup_{\epsilon>0} \limsup_{t\rightarrow\infty} \frac{1}{ t} \log  M(Y,t,\epsilon).$$

Let $(X,\Lc)$ be a  hyperbolic Riemann  surface lamination, where $X$ is a metric space endowed with a distance $\dist_X$. 
Recall that  for every $x\in X,$  a universal covering map $\phi_x$  of the leaf $L_x$  with $\phi_x(0)=x$ is  given  in \eqref{e:covering_map}.  For every  universal covering map $\psi$
of the leaf $L_x$ with $\psi(x)=x,$ there is $\theta\in\R$ such that $\psi$ is equal to  the map
 $\D\ni\xi\mapsto \phi_x(e^{i\theta}\xi),$ in other words,
those maps $\psi$  are unique up to a rotation on $\D.$
 Define the family of distances  $\Dc:=\{\dist_t\}:$
\begin{equation}\label{e:h-dist}
\dist_t(x,y):=\inf\limits_{\theta\in \R} \sup_{\xi\in \D_t} \dist_X(\phi_x(e^{i\theta}\xi),\phi_y(\xi)).
\end{equation}
The metric $\dist_t$ measures how  far two leaves
 get apart before  the hyperbolic time $t.$ It takes into account the time parametrization like in the classical
case where one measures the distance of two orbits before time $n$, by
measuring the distance at each time $i<n$. So, we  are not just concerned
 with geometric proximity.
 \begin{definition}\label{D:h-entropy}{\rm  (Dinh-Nguyen-Sibony \cite{DinhNguyenSibony14a}).} \rm
  The  {\it hyperbolic entropy} of a hyperbolic  Riemann surface $(X,\Lc),$  denoted by $h(\Lc),$ is  the entropy, computed by
\eqref{e:entropy}, of    $X$ with respect to  the family $\Dc:=\{\dist_t\}$ which is given by \eqref{e:h-dist}.
 \end{definition}
So, the value of the entropy  is  unchanged 
 under homeomorphisms  between laminations which are holomorphic along leaves. The  advantage here is that the hyperbolic time  we choose is canonical. 
 The notion of hyperbolic entropy can be extended to ($l$-dimensional) Riemannian foliations with  singularities, or more generally ($l$-dimensional) laminations with singularities, and
 a priori it is bigger than or equal to the geometric entropy of 
 Ghys, Langevin and Walczak  (see \cite[p. 584]{DinhNguyenSibony14a} for the  definition of the latter entropy and  its  relation with the  hyperbolic entropy).

\begin{theorem}\label{T:finite_entropy_regular}{\rm (Dinh-Nguyen-Sibony \cite{DinhNguyenSibony14a}).}  
Let $(X,\Lc)$ be a transversally smooth   
compact lamination  by hyperbolic  Riemann surfaces. Embed the lamination in an $\R^N$ in order to use the distance $\dist_X$ induced by a Riemannian metric on $\R^N$. 
Then,  $2\leq h(\Lc)<\infty.$ 
\end{theorem}
 
 The following proposition gives a simple criterion for the finiteness of the  hyperbolic  entropy.
We will need it  for the proof of  Theorem \ref{T:finite_entropy_regular}.
 
 \begin{proposition} \label{P:abstract_finite_entropy}
 Assume that there are positive constants $A$ and $m$ such that
for every $\epsilon>0$ small enough $X$ admits a covering by less than $A\epsilon^{-m}$ 
balls of radius $\epsilon$ for the distance $\dist_X$. Assume also that  
$$\dist_t\leq e^{ct+d}\dist_X+\varphi(t)$$
for some constants $c,d\geq 0$ and  a function $\varphi$  with $\varphi(t)\to 0$ as $t\to\infty$. Then, the entropy $h_\Dc$ is at most equal to $mc$.
\end{proposition}

 We  will apply  Proposition \ref{P:abstract_finite_entropy} in order to prove Theorem \ref{T:finite_entropy_regular}.
 So  it is  necessary 
 to  estimate the distance $\dist_t$ between leaves.
 For this  purpose, we 
 use the Beltrami equation as in the proof of the  tranverse regularity of the Poincar\'e metric in Theorem \ref{T:Poincare}.
 
 To study the   finiteness of the entropy for singular holomorphic  foliations is   a very hard matter. A satisfactory answer  is only  obtained  for complex surfaces.
 
\begin{theorem} \label{T:finite_entropy_singular} {\rm (Dinh-Nguyen-Sibony \cite{DinhNguyenSibony14b})}  Let $(X,\Lc,E)$ be a singular foliation by
Riemann surfaces  on a compact Hermitien complex  surface $X$. Assume that the singularities are linearizable and that
 the foliation is  Brody hyperbolic.
Then, its hyperbolic entropy $h(\Lc)$ is  finite.
\end{theorem}

 The proof of this theorem  is  quite delicate and requires  a  careful analysis  of the dynamics around the singularities.  
We deduce from the above theorem and   Theorem \ref{T:generic} the  following corollary. It can be applied to foliations of degree at least 2 with hyperbolic singularities.

\begin{corollary} \label{C:finite_entropy_singular}  {\rm (Dinh-Nguyen-Sibony \cite{DinhNguyenSibony14b}).}
Let $(\P^2,\Lc,E)$ be a singular foliation by
Riemann surfaces  on the complex projective plane $\P^2$ endowed with the  Fubini-Study metric. Assume that the singularities are linearizable. 
Then, the hyperbolic entropy $h(\Lc)$ of $(\P^2,\Lc,E)$ is  finite.  
\end{corollary}

Consider an  abstract   setting of  a  metric  space $(X,\dist_X)$ endowed with a  family $\Dc:=\{\dist_t\}_{t\geq 0}$ of distances. 
Let $m$ be a probability measure on $X$. 
For positive constants $\epsilon,\delta$ and $t,$ let $N_m(t,\epsilon,\delta)$ be the minimal
number of balls of radius $\epsilon$ relative to the metric
$\dist_t$ whose union has at least $m$-measure $1-\delta.$ 
{\it The (metric) entropy} of $m$ is    defined by the following formula
$$
h_\Dc(m):=\lim_{\delta\to 0}\lim_{\epsilon\to 0}\limsup_{t\to\infty}{1\over t} \log N_m(t,\epsilon,\delta). 
$$ 
We have the following general property.

\begin{lemma}
For any probability measure $m$  on $X$, we have 
$$h_\Dc(m)\leq  h_\Dc(\supp(m)),$$   where  $h_\Dc(\supp(m))$ is  computed using  formula \eqref{e:entropy}.
In particular, the metric  entropy of a probability measure is  dominated by the entropy
of the  whole  space.
\end{lemma}

As in Brin-Katok's theorem \cite{BrinKatok}, we can introduce the {\it local entropies} of $m$ at 
$x\in X$ by
$$h_\Dc^+(m,x,\epsilon):=\limsup_{t\rightarrow \infty} -{1\over t}\log m(B_t(x,\epsilon)),\qquad
h_\Dc^+(m,x):=\sup_{\epsilon> 0} h_\Dc^+(m,x,\epsilon), $$
and
$$h_\Dc^-(m,x,\epsilon):=\liminf_{t\rightarrow \infty} -{1\over t}\log m(B_t(x,\epsilon)),
\qquad
h_\Dc^-(m,x):=\sup_{\epsilon> 0} h_\Dc^-(m,x,\epsilon),
$$
where $B_t(x,\epsilon)$ denotes the ball centered at  $x$ of radius $\epsilon$ with respect to the distance $\dist_t$.

Note that
in the case of an  ergodic invariant measure associated with a continuous map on a metric compact space, the above notions of entropies coincide with the classical entropy of $m$, see Brin-Katok \cite{BrinKatok}.

Recall that in \eqref{e:h-dist} we have associated to $(X,\Lc)$ a special family of distances $\Dc:=\{\dist_t\}_{t\geq 0}$. 
Therefore, we can associate to $m$ a metric entropy and local entropies defined as above in the abstract setting. 
Recall also that a harmonic probability measure $m$ is {\it extremal} if all harmonic probability measures $m_1,m_2$ satisfying $m_1+m_2=2m$ are equal to $m$. We have the following result.

\begin{theorem} \label{th_local_entropy}{\rm (Dinh-Nguyen-Sibony \cite{DinhNguyenSibony14a}).}
Let $(X,\Lc)$ be a compact  transversally smooth  lamination by hyperbolic Riemann surfaces. 
Let $m$ be  a  harmonic probability measure.  Then, the local entropies
$h^\pm$ of $m$ are constant on leaves.  
In particular, if $m$ is  extremal,  then  $h^\pm$  are constant $m$-almost everywhere.
\end{theorem}

We  reproduce  from \cite{DinhNguyenSibony14a}  some fundamental problems concerning metric entropies for Riemann surface laminations. 
Assume, in the  first two  problems that $(X,\Lc)$ is a compact transversally smooth lamination by hyperbolic Riemann surfaces. However,   the problems can be stated in a more general setting. 

\begin{problem}\rm
Consider extremal harmonic probability measures $m$.  Is the following {\it variational principle} always true
$$h(\Lc)=\sup_m h(m)\quad  ?$$
\end{problem}

Even when this principle does not hold, it would be of interest  to consider the quantity 
$$h(\Lc)-\sup_m h(m)$$
and to explain  the role of the hyperbolic time in this number.

\begin{problem}\rm
If $m$ is as above, is the identity $h^+(m)=h^-(m)$ always true ?
\end{problem}

We think that the answer is affirmative and gives an analog of the Brin-Katok theorem.

Notice that there is a notion of entropy for harmonic measures introduced by Kaimanovich \cite{Kaimanovich}.
Consider a metric $g$ of bounded geometry on the leaves of the lamination. Then,  
we can consider the heat kernel $p(t,\cdot,\cdot)$ associated to the Laplacian
determined by this metric.  If $m$ is  a  harmonic probability measure on $X,$    Kaimanovich defines the entropy
of $m$ as
$$
h_K(m):=\int  dm(x)\Big( \lim\limits_{t\to\infty} -{1\over t}\int p(t,x,y)\log p(t,x,y)g(y)\Big).
$$ 
He  shows that the limit  exists and is  constant $m$-almost everywhere when $m$ is  extremal. 

This notion of entropy has  been extensively  studied for the universal covering of a compact Riemannian manifold, see e.g. Ledrappier \cite{Ledrappier}. 

\begin{problem}\rm
It is interesting to find relations between  Kaimanovich  entropy and   our notions of entropy. Moreover, studying these  relations   in the context of    singular holomorphic  foliations
is also an important question.
\end{problem}
In  Kaimanovich's entropy,  the  transverse  spreading is  present through the variation of the heat kernel from leaf to leaf. Therefore, a  question naturally arises whether
one can make this dependence more explicit.

Here is an open  problem  from \cite{DinhNguyenSibony14b}. 

\begin{problem}\rm
If $\Fc$ is a generic element in $ \Fc_d(\P^k)$  (with $d>1$ and $k>2$),  is   the  hyperbolic entropy of $\Fc$   finite ?
The   same question  is  asked  for a singular  holomorphic foliation  on a compact Hermitian complex manifold.
This is  a  generalization  of Theorem \ref{T:finite_entropy_singular} to higher dimensions. 
\end{problem}

Finally,  the following  problem  seems of interest.

\begin{problem}\rm
 Is  the  lower bound in  
Theorem \ref{T:finite_entropy_regular} optimal ? If yes,  study  the class of all  
 transversally smooth   
compact lamination $(X,\Lc)$  by hyperbolic  Riemann surfaces 
such that  $ h(\Lc)=2.$ 
\end{problem}

\section{Lyapunov theory for hyperbolic Riemann surface laminations} \label{S:Lyapunov}

The purpose of this section is to present some  recent results obtained in our   works \cite{NguyenVietAnh17a,NguyenVietAnh17b}.
\subsection{Brownian motion and Wiener measures}
We   start with 
 Garnett's theory   of leafwise  Brownian  motion in \cite{Garnett} (see also \cite{Candel2,CandelConlon2}).
 Our presentation follows  \cite{NguyenVietAnh17b}.
 We first recall the construction of the Wiener measure $W_0$ on the  Poincar\'e  disc $(\D,g_P).$
 Let $\Omega_0$ be  the space consisting of  all continuous  paths  $\omega:\ [0,\infty)\to  \D$  with $\omega(0)=0.$
A {\it  cylinder  set (in $\Omega_0$)} is a 
 set of the form
$$
C=C(\{t_i,B_i\}:1\leq i\leq m):=\left\lbrace \omega \in \Omega_0:\ \omega(t_i)\in B_i, \qquad 1\leq i\leq m  \right\rbrace,
$$
where   $m$ is a positive integer  and the $B_i$'s are Borel subsets of $\D,$ 
and $0< t_1<t_2<\cdots<t_m$ is a  set of increasing times.
In other words, $C$ consists of all paths $\omega\in  \Omega_0$ which can be found within $B_i$ at time $t_i.$
Let $\Ac_0$ be the  $\sigma$-algebra on $\Omega_0$  generated  by all  cylinder sets.
 For each cylinder  set  $C:=C(\{t_i,B_i\}:1\leq i\leq m)$ as  above, define
\begin{equation}\label{eq_formula_W_x_without_holonomy}
W_x(C ) :=\Big (D_{t_1}(\textbf1_{B_1}D_{t_2-t_1}(\textbf1_{B_2}\cdots\textbf1_{B_{m-1}} D_{t_m-t_{m-1}}(\textbf1_{B_m})\cdots))\Big) (x),
\end{equation}
where,  $\textbf1_{B_i}$
is the characteristic function of $B_i$ and $D_t$ is the diffusion operator
given  by  (\ref{e:diffusions}) where  $p(x,y,t)$ therein is  the heat kernel     of the Poincar\'e disc. 
 It is  well-known that $W_0$ can be   extended   to a unique probability measure on $(\Omega_0,\Ac_0).$
This  is the {\it canonical  Wiener measure}  at $0$  on the Poincar\'e disc.

Let $(X,\Lc)$ be  a hyperbolic Riemann surface lamination endowed  with the leafwise Poincar\'e metric $g_P.$
Let  $\Omega:=\Omega(X,\Lc) $  be  the space consisting of  all continuous  paths  $\omega:\ [0,\infty)\to  X$ with image fully 
contained  in a  single   leaf. This  space  is  called {\it the sample-path space} associated to  $(X,\Lc).$
  Observe that
$\Omega$  can be  thought of  as the  set of all possible paths that a 
Brownian particle, located  at $\omega(0)$  at time $t=0,$ might  follow as time  progresses.
For each $x \in  X,$
let $\Omega_x=\Omega_x(X,\Lc)$ be the  space  of all continuous
leafwise paths  starting at $x$ in $(X,\Lc),$ that is,
$$
\Omega_x:=\left\lbrace   \omega\in \Omega:\  \omega(0)=x\right\rbrace.
$$
For each $x\in X,$    the  following  mapping 
\begin{equation}\label{e:Omega_0_vs_Omega_x}
\Omega_0\ni\omega \mapsto  \phi_x\circ\omega\quad
 \text{maps}\quad  \Omega_0\quad\text{bijectively onto}\quad \Omega_x,
 \end{equation}
 where   $\phi_x:\D\to L_x$ is given in (\ref{e:covering_map}).   
Using this bijection  
we obtain   a  natural $\sigma$-algebra    $ \Ac_x$ on  the  space $\Omega_x,$ and  a natural    probability (Wiener) measure $W_x$  on $\Ac_x$  as follows:
\begin{equation}\label{e:W_x}
 \Ac_x:=\{ \phi_x\circ A:\  A\in \Ac_0\}\quad\text{and}\quad    W_x(\phi_x\circ A):=W_0(A),\qquad   A\in  \Ac_0,
\end{equation}
   where $\phi_x\circ A:= \{\phi_x\circ \omega:\ \omega\in A \}\subset \Omega_x.$

\subsection{Cocycles}
The notion of (multiplicative)  cocycles   have been   introduced 
in  \cite{NguyenVietAnh17a}  for ($l$-dimensional)  laminations. For  the sake of simplicity  we  only formulate this notion  for Riemann surface laminations in this article.  In the rest of the section   
 we make  the following convention: $\K$ denotes either the  field $\R$ or $\C.$
Moreover,  given any  integer $d\geq 1,$  $\GL(d,\K)$   denotes  the general linear group of degree $d$ over $\K$
and $\P^d(\K)$ denotes the  $\K$-projective space of dimension $d.$
\begin{definition}\label{D:cocycle} {\rm (Nguyen \cite[Definition 3.2]{NguyenVietAnh17a}). }\rm
 A {\it $\K$-valued   cocycle} (of rank $d$)   is
  a   map  
$\mathcal{A}:\ \Omega\times \R^+ \to  \GL(d,\K)      $
such that\\  
(1)  ({\it identity law})  
$\mathcal{A}(\omega,0)=\id$  for all $\omega\in\Omega ;$\\
(2) ({\it homotopy law}) if  $\omega_1,\omega_2\in \Omega_x$ and $t_1,t_2\in \R^+$ such that 
     $\omega_1(t_1)=\omega_2(t_2)$
and $\omega_1|_{[0,t_1]}$ is  homotopic  to  $\omega_2|_{[0,t_2]}$ (that is, the path $\omega_1|_{[0,t_1]}$ can be  deformed  continuously on  $L_x$ to the path  $\omega_2|_{[0,t_2]},$ 
the two endpoints of $\omega_1|_{[0,t_1]}$  being kept fixed  during the deformation), then 
$$
\mathcal{A}(\omega_1,t_1)=\mathcal{A}(\omega_2,t_2);
$$
(3) ({\it multiplicative law})    $\mathcal{A}(\omega,s+t)=\mathcal{A}(\sigma_t(\omega),s)\mathcal{A}(\omega,t)$  for all  $s,t\in \R^+$ and $\omega\in \Omega;$\\
(4) ({\it measurable law})  the {\it local expression} of $\mathcal{A}$ on each  laminated  chart is   Borel measurable.
Here,   the {\it local expression} of  $\mathcal A$ on the laminated chart 
$\Phi:  \U\to \D\times \T,$  is the map
$A:\  \D\times \D\times \T\to\GL(d,\K)$  defined  by
$$
A(y,z,t):=\mathcal A(\omega,1),
$$
where  $\omega$ is  any leafwise path  such that $\omega(0)=\Phi^{-1}(y,t),$ $\omega(1)=\Phi^{-1}(z,t)$ 
and  $\omega[0,1]$ is  contained in the simply connected  plaque   $\Phi^{-1}(\cdot,t).$   
 \end{definition}
 
 A cocycle  $\mathcal A$ on a  smooth Riemann surface lamination $(X,\Lc)$ is called   a {\it   smooth}  if,
for each laminated chart $\Phi,$ the local expression   $A$  is   smooth with respect to  $(y,z)$
 and its  partial  derivatives of any  total order 
 with respect to $(y,z)$ are jointly continuous   in $(y,z,t).$
  
 The  cocycles of rank $1$  have been  investigated  by several  authors (see, for example, Candel \cite{Candel2}, Deroin \cite{Deroin}, etc.).   
 The  holonomy cocycle (or equivalently the normal  derivative cocycle)  of the  regular part of a singular holomorphic  foliation  by hyperbolic Riemann  surfaces $(X,\Lc,E)$ with $\dim_\C X=n$ is
 a typical   example of  $\C$-valued cocycles of rank $n-1.$   These   cocycles   capture  the topological aspect of the considered  foliations. Moreover, 
we can produce  new  cocycles from the old ones by performing some basic operations such as the wedge  product and the tensor product
(see \cite[Section 3.1]{NguyenVietAnh17a}). 

\subsection{Oseledec multiplicative ergodic theorem}

Now  we are in the position to state the Oseledec multiplicative  ergodic theorem for  hyperbolic Riemann surface laminations.
 \begin{theorem} \label{T:VA_general}  {\rm (Nguyen \cite[Theorem 3.11]{NguyenVietAnh17a}).}
 Let $(X,\Lc)$ be  a   hyperbolic Riemann surface lamination. 
Let  $\mu$ be a  harmonic measure which is  also  ergodic.
Consider  a    cocycle
$\mathcal{A}:\ \Omega\times \R^+ \to  \GL(d,\K)    .  $   Assume that the following integrability condition  is  satisfied for some real number $t_0>0:$  
\begin{equation}\label{e:integrability}
\int_{x\in X}\big(\int_{\Omega_x} \sup_{t\in[0,t_0]}\log^+\|\mathcal A(\omega,t)\|dW_x(\omega)\big)d\mu(x)<\infty,
\end{equation}
where $\log^+:=\max(0,\log).$
Then  there exist  a leafwise saturated  Borel  set $Y\subset X$ of  total $\mu$-measure  and a number $m\in\N$  together with $m$ integers  $d_1,\ldots,d_m\in \N$  such that
the following properties hold:
\begin{itemize}
\item[(i)] For   each $x\in Y$  
 there   exists a  decomposition of $\K^d$  as  a direct sum of $\K$-linear subspaces 
$$\K^d=\oplus_{i=1}^m H_i(x),
$$
 such that $\dim H_i(x)=d_i$ and  $\mathcal{A}(\omega, t) H_i(x)= H_i(\omega(t))$ for all $\omega\in  \Omega_x$ and $t\in \R^+.$   
Moreover,  $x\mapsto  H_i(x)$ is   a  measurable map from $  Y $ into the Grassmannian of $\K^d.$
For each $1\leq i\leq m$ and each $x\in Y,$ let $V_i(x):=\oplus_{j=i}^m H_j(x).$  Set $V_{m+1}(x)\equiv \{0\}.$
\item[
(ii)]  There   are real numbers 
$$\chi_m<\chi_{m-1}<\cdots
<\chi_2<\chi_1,$$
  and   for  each $x\in Y,$ there is a set $F_x\subset \Omega_x$ of total $W_x$-measure such that for every $1\leq i\leq m$ and  every  $v\in V_i(x)\setminus V_{i+1}(x)$
 and every  $\omega\in F_x,$
\begin{equation}
\label{e:Lyapunov}
\lim\limits_{t\to \infty, t\in \R^+} {1\over  t}  \log {\| \mathcal{A}(\omega,t)v   \|\over  \| v\|}  =\chi_i.    
\end{equation}
Moreover, 
\begin{equation}
\label{e:Lyapunov_max}
\lim\limits_{t\to \infty, t\in \R^+} {1\over  t}  \log {\| \mathcal{A}(\omega,t)  \|}  =\chi_1    
\end{equation}
 for  each $x\in Y$  and for   every  $\omega\in F_x.$ 
\end{itemize}
Here    $\|\cdot\|$  denotes the standard   Euclidean norm of $\K^d.$  
 \end{theorem}  
The  above  result  is  the counterpart, in the context of hyperbolic Riemann surface laminations,
of the classical  Oseledec   multiplicative ergodic theorem for maps   (see \cite{KatokHasselblatt,Oseledec}).
In fact, Theorem 3.11 in  \cite{NguyenVietAnh17a} is  much more general than  Theorem \ref{T:VA_general}.
Indeed, the former is formulated for $l$-dimensional  laminations and for  leafwise Riemannian metrics which satisfy some
reasonable geometric   conditions.
  
 Assertion  (i) above tells us that the  Oseledec  decomposition exists for all points $x$ in  a leafwise saturated Borel  set
of  total $\mu$-measure and that this  decomposition is holonomy invariant.
Observe that  the Oseledec  decomposition  in  (i)  depends only on $x\in Y,$ in particular, it 
does not depend on paths $\omega\in\Omega_x.$

The decreasing  sequence  of  subspaces  of $\K^d$ given by assertion (i):
$$
\{0\}\equiv V_{m+1}(x)\subset V_m(x)\subset \cdots \subset V_1(x)=\K^d
$$
is  called the {\it Lyapunov filtration} associated to $\mathcal A$  at a given point $x\in Y.$ 
 
 The   numbers  $\chi_m<\chi_{m-1}<\cdots
<\chi_2<\chi_1$ given by  assertion (ii) above are called  the {\it  Lyapunov exponents} of the cocycle $\mathcal{A}$ with respect to the harmonic measure $\mu.$
  It follows   from formulas (\ref{e:Lyapunov}) and (\ref{e:Lyapunov_max}) above  that these characteristic numbers measure heuristically the  expansion rate of
  $\mathcal A$ along different vector-directions $v$ and  along leafwise Brownian trajectories.
In other words,  the stochastic formulas (\ref{e:Lyapunov})-(\ref{e:Lyapunov_max})  highlight the  dynamical  character  of the Lyapunov exponents.

\subsection{Applications to compact smooth laminations and  compact singular foliations}\label{SS:}

Let $\mathcal{A}:\ \Omega\times \R^+ \to  \GL(d,\K)$ be a smooth cocycle  defined on a  smooth hyperbolic  Riemann surface lamination $(X,\Lc).$ 
Observe that  the map $\mathcal{A}^{*-1}:\ \Omega\times \R^+ \to  \GL(d,\K),$ defined by    
  $\mathcal{A}^{*-1}(\omega,t):= \big (\mathcal A(\omega,t)\big)^{*-1},$ is  also a cocycle,
where  $A^*$ (resp. $A^{-1}$) denotes as  usual the transpose (resp. the inverse)  of a square  matrix $A.$ 

We define  two functions  $\bar\delta(\mathcal A),\ \underline\delta(\mathcal A):\ X\to\R$ as well  as four quantities  $\bar\chi_{\max}( \mathcal A ),$
$ \underline\chi_{\max}( \mathcal A ),$    $ \bar\chi_{\min}( \mathcal A ),$
$\underline\chi_{\min}( \mathcal A )$
   as  follows.
  Fix  a point $x\in  X,$ an element $u\in\K^d\setminus \{0\}$  and     a  simply connected plaque $K$   of $(X,\Lc)$ passing through
$x.$
Consider   the  function  $f_{u,x}:\  K\to \R$ defined by
\begin{equation}\label{eq_function_f}
f_{u,x}(y):= \log {\| \mathcal A(\omega,1)u \|\over  \| u\|} ,\qquad  y\in K,\ u\in\K^d\setminus\{0\},
\end{equation}
where  $\omega\in \Omega$ is any path  such that $\omega(0)=x,$ $\omega(1)=y$ 
and that $\omega[0,1]$ is  contained in $K.$  Then define
\begin{equation}\label{eq_formulas_delta}
\bar \delta(\mathcal A)(x):=\sup_{u\in \K^d:\ \|u\|=1} (\Delta f_{u,x})(x)\ \ \text{and}\  \
 \underline \delta(\mathcal A)(x):=\inf_{u\in \K^d:\ \|u\|=1} (\Delta f_{u,x})(x),\\
\end{equation}
where $\Delta$ 
 is, as  usual,   the  Laplacian   
on the leaf $L_x$ induced  by the leafwise Poincar\'e metric  $g_P$ on  $(X,\Lc)  $ (see formula \eqref{e:Laplacian} for $\beta:=g_P$).
We also define
\begin{equation}\label{eq_formulas_chi}
\begin{split}
\bar\chi_{\max}=\bar\chi_{\max} (\mathcal A)&:=\int_X \bar\delta(\mathcal A)  (x)   d\mu(x),\\
\underline\chi_{\max}=\underline\chi_{\max} (\mathcal A)&:=\int_X \underline\delta(\mathcal A)  (x)   d\mu(x);\\
\underline\chi_{\min}=\underline\chi_{\min}( \mathcal A )&:=-    \bar\chi_{\max}(\mathcal A^{*-1}) ,\\
\bar\chi_{\min}=\bar\chi_{\min}( \mathcal A )&:=-   \underline\chi_{\max}(\mathcal A^{*-1}) .
\end{split}
\end{equation}
Note that our functions  $\bar\delta,$ $\underline\delta$ are  the multi-dimensional generalizations  of the operator $\delta$ introduced by Candel \cite{Candel2}.


We are in the position to state  effective  integral estimates on the Lyapunov  exponents. 

 \begin{theorem}  \label{T:VA_smooth}  {\rm (Nguyen \cite[Theorem 3.12]{NguyenVietAnh17a}).}
 Let $(X,\Lc)$ be   a      compact  smooth  lamination by hyperbolic Riemann surfaces. 
 Let $\mu$ be a   harmonic probability measure which is  ergodic.
Let
$\mathcal{A}:\ \Omega\times \R^+ \to  \GL(d,\K)      $ be  a smooth cocycle.  
 Let    
$$\chi_m<\chi_{m-1}<\cdots
<\chi_2<\chi_1$$
be the Lyapunov exponents of the cocycle $\mathcal A$ with respect  to $\mu,$ given by Theorem \ref{T:VA_general}.
Then  the following  inequalities hold $$
\underline\chi_{\max}\leq \chi_1\leq \bar\chi_{\max} \quad\text{and}\quad
\underline\chi_{\min} \leq \chi_m\leq \bar\chi_{\min} .$$ 
\end{theorem} 
 
 This theorem generalizes some  results of Candel \cite{Candel} and Deroin  \cite{Deroin}  who  treat the case $d=1.$
 Under the assumption of   Theorem  \ref{T:VA_smooth}, the integrability  condition  \eqref{e:integrability}
   follows from some  well-known  estimates of the heat kernels of the Poincar\'e disc and  the fact that  the lamination is  
 compact and is  without singularities. In fact, we improve  the method of Candel in \cite{Candel2}.

 The  holonomy cocycle (or equivalently, the normal derivative cocycle)  of a foliation  
is  a very important object  which  encodes
 dynamical  as well as  geometric and analytic informations    of  the  foliation.
 Exploring  this  object  allows  us  to understand  more  about the  foliation itself.
On  the other hand, recall  that the main examples of holomorphic foliations by curves 
are  those  in the complex projective space $\P^k$ of arbitrary
dimension (in which case there are always singularities) or in algebraic manifolds,
and that the typical  properties of a generic  foliation of a given degree $d>1$   in $\P^k$
are  described by Theorem \ref{T:generic}.
Therefore, 
  the following fundamental  question  arises  naturally: 

\smallskip

\noindent {\bf Question. }{\it Can one  define  the  Lyapunov exponents of 
an ergodic harmonic  measure $\mu$ on a  compact singular holomorphic hyperbolic  foliation $\Fc=(X,\Lc,E)$ ?}

\smallskip

By  Proposition \ref{P:harmonic_currents_vs_measures},  this  question can be rephrased  for  directed harmonic  currents on the foliation.
We have recently obtained  the  following  affirmative  answer  to this  question for generic  foliations in dimension two.  
 
\begin{theorem}\label{T:VA_singular}  {\rm  (Nguyen \cite[Theorem 1.1]{NguyenVietAnh18}).}
Let $\Fc=(X,\Lc,E)$    be  
a   holomorphic Brody hyperbolic foliation   with     hyperbolic  singularities $E$ in a   Hermitian compact complex projective  surface $X.$
Let $\mathcal A$ be the holonomy cocycle of the foliation. Let $T$ be  a positive harmonic  current  directed by $\Fc$  which does not give mass to any  invariant  analytic curve.
Consider   the  corresponding harmonic  measure
 $\mu:= T\wedge  g_P,$  where $g_P$ is as usual the leafwise  Poincar\'e metric. 
  Then the integrability condition \eqref{e:integrability} is  satisfied for all $t_0>0.$
\end{theorem}
Here   is an immediate  consequence of this theorem.

\begin{corollary}\label{C:extremal}
 Under the hypotheses and notation of Theorem \ref{T:VA_singular}, assume in addition that
the measure $\mu$ is  ergodic.    Then  $T$ admits
the (unique) Lyapunov exponent $\lambda(T)$  given by the  formula
$$
\lambda(T):= \int_X\big(\int_{\Omega_x} \log \|\mathcal{A}(\omega,1)\| dW_x(\omega)\big)d\mu(x).  
$$
Moreover,
  for  $\mu$-almost  every $x\in X,$  we have
 $$
 \lim\limits_{t\to \infty} {1\over  t} \log  \| \mathcal A(\omega,t) \|=\lambda(T)  
 $$
  for    almost every path  $\omega\in\Omega$ with respect to the Wiener measure  at $x$
  which lives on the leaf passing through $x.$
\end{corollary}
 
 Consider     
a   singular  foliation by curves  $\Fc=(\P^2,\Lc,E)$ on the complex projective plane $\P^2$
 such that 
all the singularities of $\Fc$   are hyperbolic and  that  $\Fc$ has no invariant algebraic  curve.
 By  Remark  \ref{R:Brody_sufficiency}  we know 
that 
 $\Fc$ is   Brody  hyperbolic.  Moreover,
 Theorem \ref{T:FornaessSibony} tells us  that the harmonic current $T$ is unique
up to a multiplicative constant. So the convex cone of all harmonic  currents of $\Fc$ is just  a real half-line,
and hence all these  currents  are extremal. Therefore, by Proposition \ref{P:harmonic_currents_vs_measures}  the measure $T\wedge g_P$ is ergodic.  
Consequently, Corollary  \ref{C:extremal} applies and   gives us the following result.
It can be   applied  to  every generic  foliation in $\P^2$ with a given degree $d>1.$
  
\begin{corollary}\label{cor_th_main}
 Let $\Fc=(\P^2,\Lc,E)$    be  
a   singular  foliation by curves on the complex projective plane $\P^2.$ Assume that
all the singularities  are hyperbolic and  that  $\Fc$ has no invariant algebraic  curve. 
Let $T$ be  the unique  harmonic   current  tangent to $\Fc$  such that  $\mu:=T\wedge g_P$ is a probability measure.
  Let $\mathcal A,$  
be as in the  statement of  Theorem \ref{T:VA_singular}. Then the conclusion of this theorem
as  well as  that  of   Corollary \ref{C:extremal} hold.
In particular, $\Fc$ admits  a unique  Lyapunov  exponent.
\end{corollary}  
  
The novelty of  the last corollary is  that the   (unique)  Lyapunov  exponent of such a foliation $\Fc$ is  intrinsic and  canonical.

The proof of Theorem \ref{T:VA_singular}  consists of two steps. Let $g_X$ be  a Hermitian  metric on $X$
 and let $\dist$ denote the  distance  on $X$ induced by $g_X.$ 
In the  first step  we show that    Theorem \ref{T:VA_singular}  follows  from
the new integrability condition  \eqref{e:necessary_integrability}.
\begin{equation}\label{e:necessary_integrability}
 \text{(new integrability  condition):}\qquad \int_X| \log \dist(x,E)| \cdot (T\wedge g_P)(x)<\infty.
 \end{equation}
 This new condition has  the advantage over  the old one \eqref{e:integrability}, since  the former does not involve the somewhat complicating  Wiener measures, and hence  it is
 easier to handle than the latter.

For  this purpose  we  study  the behavior of the holonomy cocycle
near the  singularities  with respect to  the leafwise Poincar\'e metric.  
Roughly speaking, this step
 quantifies  the  expansion speed of the hololomy  cocycle in terms
 of the ambient metric $g_X$ when  one travels  along   unit-speed geodesic rays.
 One of the  main ingredients  is a  detailed  analysis of  the  behaviour  of the leafwise Poincar\'e metric near hyperbolic  singularities
 which  has been  carried  out  in  \cite{DinhNguyenSibony12,DinhNguyenSibony14a,DinhNguyenSibony14b}.

  The  second main step is  then devoted to  the proof of  inequality \eqref{e:necessary_integrability}.
  The main difficulty is that  known  estimates (see, for example,  \cite{DinhNguyenSibony12})
 on the behavior of $T$
near  linearizable  singularities,    only  give a  weaker  inequality
  \begin{equation} \label{e:known_estimate}
\int_X | \log \dist(x,E)|^{1-\delta}\cdot( T\wedge g_P)(x)<\infty, \qquad\forall \delta>0.
  \end{equation}
  So \eqref{e:necessary_integrability} is  the limiting case of  \eqref{e:known_estimate}.
  The proof of \eqref{e:known_estimate}  relies on  the finiteness of the Lelong number of $T$  at
  every   point which has been   established in Proposition \ref{P:Skoda}. 
 Recall that Theorem  \ref{T:Lelong}   sharpens the last estimate
by showing that the Lelong number  of $T$ vanishes at every   hyperbolic singular point $x\in E.$  
Nevertheless, even this better estimate does not suffice to prove \eqref{e:necessary_integrability}.

  The new idea  in \cite{NguyenVietAnh18} is  that  we use  a cohomological argument
  which  exploits  fully  the assumption that $X$ is  projective.
  This  assumption  imposes a stronger mass-clustering  condition on harmonic  currents  than \eqref{e:known_estimate}.

The  condition  of Brody hyperbolicity  seems to be indispensable  for the   integrability of the  holonomy cocycle. Indeed,
 Hussenot  \cite[Theorem A]{Hussenot} finds out    the  following  remarkable property
 for  a class of  Ricatti foliations $\Fc$ on $\P^2.$
For every $x\in\P^2$ outside invariant curves of every foliation in this class, it holds that
   $$
 \limsup\limits_{t\to \infty} {1\over  t} \log  \| \mathcal A(\omega,t) \|= \infty 
 $$ for almost every path $\omega\in\Omega_x$ with respect to the Wiener measure at $x$ which lives
 on the leaf passing through $x.$ By     Theorem  \ref{T:generic}, these foliations  are hyperbolic since all their  singular points    have nondegenerate linear part.
 Nevertheless, neither of them is Brody hyperbolic because  they all contain integral curves which are some images of $\P^1$  (see  Remark  \ref{R:Brody_sufficiency}).

  \begin{problem}\rm {\rm (see also  \cite{DeroinKleptsyn,Hussenot}).}
   Is  the Lyapunov exponent of a generic  foliation with a given degree $d>1$ in $\P^2$  positive/negative/zero ?
  \end{problem}

  \begin{problem}\rm 
Does  Theorem \ref{T:VA_singular}   still hold if 
the ambient compact projective manifold $X$ is of dimension $>2$ ? 
  \end{problem}

\subsection{Geometric   characterization of Lyapunov exponents}

To find  a  geometric  interpretation  of these  characteristic  quantities,
 our idea  consists    in replacing the Brownian trajectories by the more appealing objects, namely, the {\it     unit-speed   geodesic rays.}  These paths   are  parameterized  by  their  length (with respect to the leafwise Poincar\'e metric).   Therefore, we characterize 
 the Lyapunov exponents  in terms of  the  expansion rates of $\mathcal A$ along the geodesic rays.
 %
 
 Let $(X,\Lc)$ be  a   hyperbolic Riemann surface lamination.
 Recall from (\ref{e:covering_map}) that
 $(\phi_x)_{x\in X}$ is  a given family of  universal covering maps $\phi_x:\ \D\to L_x$  with $\phi_x(0)=x.$
 For every $x\in X,$  
  the  set of all unit-speed geodesic rays $\omega:\ [0,\infty)\to L_x$ starting at $x$ (that is,    $\omega(0)=x$),     can  be described   by  the family  $(\gamma_{x,\theta})_ {\theta\in [0,1)},$ where
 \begin{equation}\label{eq_geodesics}
 \gamma_{x,\theta}(R):=  \phi_x(e^{2\pi i\theta}r_R),\qquad  R\in\R^+,
 \end{equation}
 and  $r_R$  is uniquely determined  by the  equation $r_R\D=\D_R $  (see \eqref{e:radii_conversion}).
 The path $\gamma_{x,\theta}$ is called the {\it  unit-speed geodesic ray}   at $x$ with the leaf-direction $\theta.$
 Unless  otherwise  specified, the {\it space of leaf-directions} $[0,1)$ is  endowed with the Lebesgue measure. 
 This space   is  visibly identified, via the map $[0,1)\ni\theta\mapsto e^{2\pi i\theta},$ with the unit circle $\partial \D$ endowed with the normalized rotation measure.

 We  introduce the following   notions of expansion  rates for  cocycles.
 
  \begin{definition}\label{D:expansion_rate}\rm
 Let $\mathcal A$ be a $\K$-valued   cocycle and  $R>0$  a  time.

  The {\it expansion rate}  of $\mathcal A$  at a  point $x\in X$ in the leaf-direction $\theta$  at time  $R$ along the vector $v\in \K^d\setminus\{0\}$ 
 is the number
 $$\Ec(x,\theta,v,R):={1\over  R}  \log {\| \mathcal A(\gamma_{x,\theta},R)v   \|\over  \| v\|}.$$ 
 
 The {\it  expansion rate}  of $\mathcal A$   at a  point $x\in X$ in the leaf-direction $\theta$ 
  at time  $R$
 is
 \begin{equation*}\begin{split}
 \Ec(x,\theta,R):= \sup\limits_{v \in \K^d\setminus\{0\}} \Ec(x,\theta,v,R)&=\sup_{v \in \K^d\setminus\{0\} } {1\over  R}  \log {\| \mathcal A(\gamma_{x,\theta},R)v   \|\over  \| v\|}\\
 &=
 {1\over  R}  \log {\| \mathcal A(\gamma_{x,\theta},R)  \|}.
 \end{split}
 \end{equation*}

 Given a $\K$-vector subspace   $\{0\}\not=H\subset \K^d,$ the {\it expansion rate}  of $\mathcal A$  at a  point $x\in X$ at time  $R$
 along the vector space  $H$  
 is the interval $\Ec(x,H,R):=[a,b],$ where
 $$
 a:=  \inf_{v \in H\setminus \{0\}} \int_0^1   \Big ({1\over  R}  \log {\| \mathcal A(\gamma_{x,\theta},R)v   \|\over  \| v\|}\Big) d\theta\ \text{and}\
 b:=  \sup_{v \in H\setminus \{0\}} \int_0^1   \Big ({1\over  R}  \log {\| \mathcal A(\gamma_{x,\theta},R)v   \|\over  \| v\|}\Big) d\theta.
 $$
 \end{definition} 
  Notice  that    $ 
 \Ec(x,\theta,v,R) $ (resp. $
 \Ec(x,\theta,R)$)  expresses  geometrically the  expansion rate (resp. the maximal expansion rate) of the  cocycle 
when one travels along the unit-speed geodesic ray $\gamma_{x,\theta}$ up to time $R.$  
  On the other hand, $ 
 \Ec(x,H,R) $   represents the  smallest closed interval which  contains all numbers
$$\int_0^1   \Big ({1\over  R}  \log {\| \mathcal A(\gamma_{x,\theta},R)v   \|\over  \| v\|}\Big) d\theta,
$$
where  $v$ ranges over $H\setminus \{0\}.$
Note that  the above integral is the average of  the expansion rate  of the  cocycle 
when one travels along the unit-speed geodesic rays along the vector-direction $v\in H$   from $x$  to the  Poincar\'e circle  with radius $R$  and  center $x$ spanned  on $L_x.$ 

We say  that a sequence of  intervals $[a(R),b(R)]\subset \R$ indexed by $R\in\R^+$ converges to  a number $\chi\in \R$ and write $\lim_{R\to\infty} [a(R),b(R)]=\chi,$ if   
  $\lim_{R\to\infty} a(R)= \lim_{R\to\infty} b(R)=\chi.$
  
  Now  we are able  to state the main result of this  subsection.
  
  \begin{theorem} \label{T:Geometric_Lyapunov} {\rm (Nguyen \cite{NguyenVietAnh17b}).}
   Let $(X,\Lc)$ be  a compact smooth hyperbolic Riemann surface lamination and $T$  a directed  positive harmonic current
   which is also extremal. 
 Let 
  $\mu:=T\wedge g_P$ be the (positive finite Borel)  measure associated to $T.$  
Consider  a   smooth  cocycle
$\mathcal{A}:\ \Omega\times \R^+ \to  \GL(d,\K)  .    $ 
Then there is a leafwise saturated Borel set $Y$  of total $\mu$-measure 
which satisfies the conclusion of  Theorem \ref{T:VA_general}
and  the following additional geometric  properties: 
 \begin{itemize}
\item[(i)] For each $1\leq i\leq m$ and  for each $x\in Y,$ there is a  set $G_x\subset [0,1)$ of total Lebesgue measure 
 such that     for each  $v\in V_i(x)\setminus V_{i+1}(x),$ 
\begin{equation*}
 \lim_{R\to\infty}\Ec(x,\theta,v,R)= \chi_i,\qquad \theta\in G_x.
\end{equation*}
Moreover, the maximal Lyapunov exponent $\chi_1$ satisfies
 \begin{equation*}
 \lim_{R\to\infty}\Ec(x,\theta,R)= \chi_1,\qquad \theta\in G_x.
\end{equation*}  
\item[(ii)]  For  each  $1\leq i\leq m$ and each $x\in Y,$ 
\begin{equation*}
\lim_{R\to\infty}\Ec(x,H_i(x),R)=\chi_i.
\end{equation*} 
  \end{itemize}
Here
$\K^d=\oplus_{i=1}^m H_i(x),$ $x\in Y,$ is  the Oseledec  decomposition given by Theorem \ref{T:VA_general}   and      $\chi_m<\chi_{m-1}<\cdots
<\chi_2<\chi_1$ are the  corresponding     Lyapunov exponents.
\end{theorem}

 
 Theorem  \ref{T:Geometric_Lyapunov} gives a geometric meaning to  
 the stochastic formulas (\ref{e:Lyapunov})--(\ref{e:Lyapunov_max}).

 Let  $\Fc=(M,\Lc,E)$ be  a transversally smooth (resp.   holomorphic)    singular foliation  by  Riemann surfaces 
in a Riemannian manifold (resp. Hermitian complex  manifold) $M.$
Consider  a leafwise saturated, compact   set $X\subset M\setminus E$  whose  leaves are all hyperbolic.  
So the  restriction of the foliation  $(M\setminus E,\Lc)$ to $X$ gives an inherited compact
  smooth hyperbolic Riemann lamination $(X,\Lc).$ Moreover,
the   holonomy  cocycle of $(M\setminus E,\Lc)$ induces, by restriction, an inherited smooth cocycle on $(X,\Lc).$
Hence, Theorem \ref{T:Geometric_Lyapunov} applies  to  the latter cocycle.   
Recall  that
    a {\it minimal set} of $\Fc$ is a  leafwise saturated  subset  of $M\setminus  E$  which is  also a closed  subset
of $M$  and  which  contains  no proper subset with these properties. 
In particular, when $(M,\Lc,E)$ is a singular  holomorphic     foliation  on a compact Hermitian complex  manifold $M$ of dimension $n,$   
 the last theorem
applies to the induced  holonomy cocycle of rank $n-1$ associated with every  minimal set  $X$ whose leaves are all  hyperbolic.

The proof of Theorem  \ref{T:Geometric_Lyapunov} (i)  relies on the theory of  Brownian trajectories  on hyperbolic spaces.
More  concretely, 
some quantitative results on  the boundary behavior of   Brownian trajectories  
by  Lyons \cite{Lyons} and Cranston \cite{Cranston} and on  
the shadow of   Brownian trajectories   by geodesic rays  
are  our main  ingredients.
This 
  allows us to replace a  Brownian  trajectory  by a unit-speed geodesic ray with uniformly  distributed leaf-direction.
Hence,   Part (i) of  Theorem  \ref{T:Geometric_Lyapunov}  will follow from  Theorem \ref{T:VA_smooth}. 

To establish Part (ii) of Theorem \ref{T:Geometric_Lyapunov} we  need  two steps.  
In the  first step we  adapt   to our context 
the  so-called {\it Ledrappier  type  characterization  of Lyapunov  spectrum} which  was introduced in  \cite{NguyenVietAnh17a}.  
 This  allows us to show that
 a  similar version of  Part (ii)  of Theorem \ref{T:Geometric_Lyapunov}  holds  when  the expansion  rates in terms of  geodesic rays
are replaced by some  heat diffusions associated  with the  cocycle.
The  second step shows that  the above  heat diffusions can be  approximated by  the expansion  rates. 
To  this end  we    establish a new  geometric estimate on the heat diffusions
which   relies on  the  proof of the geometric Birkhoff ergodic theorem  (Theorem \ref{T:geometric_ergodic}).

\begin{problem} \rm Is  Theorem \ref{T:Geometric_Lyapunov} still true if
   $(X,\Lc)$  is  the whole regular part
of a  singular holomorphic  foliation $\Fc$ by hyperbolic Riemann surfaces on a compact complex manifold $M$ and $\mathcal A$ is  the  holonomy cocycle  ? 
We  can begin  with  the case  where $M$ is a  surface.
 \end{problem}

\end{document}